\documentclass{amsart}
\usepackage{amssymb}
\newtheorem{thm}{Theorem}[section]
\newtheorem{cor}[thm]{Corollary}
\newtheorem{lem}[thm]{Lemma}
\newtheorem{prop}[thm]{Proposition}
\theoremstyle{definition}
\newtheorem{defn}[thm]{Definition}
\theoremstyle{remark}
\newtheorem{rem}[thm]{Remark}
\numberwithin{equation}{section}
\newcommand{\To}{\longrightarrow}

\newcommand{\inv}{^{-1}}
\newcommand{\C}{\mathbb C}
\newcommand{\Z}{\mathbb Z}
\newcommand{\R}{\mathbb R}
\newcommand{\N}{\mathbb N}

\newcommand{\X}{\mathcal X}
\newcommand{\Y}{\mathcal Y}

\newcommand{\x}{\times}

\newcommand{\la}{\lambda}
\newcommand{\al}{\alpha}
\newcommand{\si}{\sigma}

\newcommand{\Id}{\operatorname{Id}}
\newcommand{\tr}{\operatorname{tr}}
\newcommand{\Ad}{\operatorname{Ad}}

\begin{document}

\title[Totally positive Toeplitz matrices and quantum cohomology]
{Totally positive Toeplitz matrices and quantum cohomology of
partial flag varieties}
\author{Konstanze Rietsch}%
\address{Mathematical Institute, University of Oxford}%
\email{rietsch@maths.ox.ac.uk}%


\thanks{The author was funded for most of this work
by EPSRC grant GR/M09506/01, and is currently supported by the Glasstone
Foundation
}


\begin{abstract}
We show that the set of totally positive unipotent
lower-triangular Toeplitz matrices in $GL_n$ form a real
semi-algebraic cell of dimension $n-1$. Furthermore we prove a
natural cell decomposition for its closure. The proof uses
properties of the quantum cohomology rings of the partial flag
varieties of $GL_n(\C)$ relying in particular on the positivity of
the structure constants, which are enumerative 
Gromov--Witten invariants. We also give
a characterization of total positivity for Toeplitz matrices in
terms of the (quantum) Schubert classes. This work builds on some
results of Dale Peterson's which we explain with proofs in the
type $A$ case.
\end{abstract}
\maketitle
\section{Introduction}
A matrix is called totally nonnegative if all of its minors are
nonnegative. Totally nonnegative infinite Toeplitz matrices were
studied first in the 1950's. They are characterized in the following
theorem conjectured by Schoenberg and
proved
by Edrei.

\begin{thm}\cite{Edr:ToeplMat} The $\infty\x\infty$--Toeplitz matrix
\begin{equation*}
A=
\begin{pmatrix}
 1        &       &       &      &         &      &       & \\
\overset{\ } a_1       &1     &       &       &         &      &      & \\
a_2       & a_1  &1      &      &         &      &       & \\
  \vdots  &  a_2   &  a_1  &\ddots &         &      &      & \\
 a_d      &       & \ddots  &\ddots &  \ddots &       &     & \\
 a_{d+1}    &\ddots &       & \ddots     &     a_1 & 1     &      & \\
 \vdots         &   \ddots    &\ddots    &      &      a_2   & a_1
&\ddots     & \\
          &       &\ddots       &\ddots &        &     \ddots  &\ddots &\ddots
\end{pmatrix}
\end{equation*}
is totally nonnegative precisely if its generating function is of
the form,
\begin{equation*}
1+ a_1 t + a_2 t^2 + \dotsc =\exp{(t
\alpha)}\prod_{i\in\N}\frac{(1 +\beta_i t)}{(1- \gamma_i t)},
\end{equation*}
where $\alpha\in \R_{\ge 0}$ and $ \beta_1\ge\beta_2\ge\cdots \ge
0, \ \gamma_1\ge\gamma_2\ge\cdots\ge 0 $ with
$\sum\beta_i+\sum\gamma_i<\infty$.
\end{thm}
This beautiful result has been reproved many times, see
\cite{Ok:Sinfty} for an overview. It may be thought of as giving a
parameterization of the totally nonnegative Toeplitz matrices by
\begin{equation*}
\{(\alpha;(\tilde\beta_i)_i,(\tilde\gamma_i)_i)\in \R_{\ge
0}\x\R_{\ge 0}^{\N}\x\R_{\ge 0}^\N\ |\ \sum_{i\in \N}
i(\tilde\beta_i + \tilde\gamma_i)<\infty\ \},
\end{equation*}
where $\tilde \beta_i=\beta_{i}-\beta_{i+1}$ and $\tilde
\gamma_i=\gamma_{i}-\gamma_{i+1}$.

Now let $U^-$ denote the lower triangular unipotent $n\x
n$--matrices. One aim of this paper is to parameterize the set of
totally nonnegative matrices in
\begin{equation*}
X:=\left \{ x\in U^-\left. |\  x=
\begin{pmatrix}
 1        &       &       &      &         &              \\
a_1       &1      &       &       &         &             \\
a_2       &  a_1  &1      &      &         &            \\
a_3       & \ddots      &  \ddots  &\ddots &         &           \\
\vdots    & \ddots       &   \ddots    &a_1 &  1 &            \\
 a_{n-1}  &\dots       & a_3      &a_2 &   a_1    & 1
 \end{pmatrix}
\right .\right\}
\end{equation*}
by $n-1$ nonnegative parameters. Let $\Delta_{n-i}  (x)$ be the
lower left hand corner $i\x i$ minor of $x\in X$. Explicitly, we
will prove the following statement.

\begin{prop}\label{t:2}
Let $X_{\ge 0}$ denote the set of totally nonnegative matrices in
$X$.  Then the map
\begin{equation*}
\Delta_{\ge 0} := (\Delta_1,\dotsc, \Delta_{n-1}):\ X_{\ge
0}\longrightarrow \R_{\ge 0}^{n-1}
\end{equation*}
is a homeomorphism.
\end{prop}
Note that $\Delta:X\to\C^{n-1}$ is a ramified cover and the
nonnegativity of the values of the $\Delta_i$ is by no means
sufficient for an element $u$ to be totally nonnegative. The
statement is rather that for prescribed nonnegative values of
$\Delta_1,\dotsc,\Delta_{n-1}$, among all matrices with these
fixed values there is precisely one which is totally nonnegative.

The proof of this result involves relating total positivity for
these $n\x n$ Toeplitz matrices to properties of quantum
cohomology rings of partial flag varieties, via Dale Peterson's
realization of these as coordinate rings of certain remarkable
subvarieties of the flag variety. We show that the Schubert basis
of the quantum cohomology ring plays a similar role for these
matrices with regard to positivity as does the (classical limit
of the) dual canonical basis for the whole of $U^-$ in the work
of Lusztig \cite{Lus:TotPos94}. This is the content of
Theorem~\ref{t:main}, which is the main result of this paper. The
above parameterization of $X_{\ge 0}$ comes as a corollary.

\subsection{Overview of the paper}

The first part of the paper is taken up with introducing the
machinery we will need to prove our results. We set out by
recalling background on the quantum cohomology rings of full and
partial flag varieties, especially work of Astashkevich, Sadov,
Kim, and Ciocan-Fontanine, as well as Fomin, Gelfand and
Postnikov.

Their work is then used in Section~\ref{s:Peterson} to explain
Peterson's result identifying these rings as coordinate rings of
affine strata of a certain remarkable subvariety $\Y$ of the flag
variety.
The variety $X$ of Toeplitz matrices enters the picture when the
Peterson variety $\Y$ is viewed from the opposite angle
($U^-$-orbits rather than $U^+$-orbits). We recall the Bruhat
decomposition of the variety of Toeplitz matrices. Each stratum
$X_P$ has in its coordinate ring the quantum cohomology ring of a
partial flag variety $G/P$ with its Schubert basis and quantum
parameters.

In Section~\ref{s:Deltas} we recall Kostant's formula for the
quantum parameters as functions on $X_B$ in terms of the minors
$\Delta_i$ and generalize it to the partial flag variety case.

After some motivation from total positivity the main results are
stated in Section~\ref{s:main}. The Theorem~\ref{t:main} has
three parts.  Firstly, the set of points in $X_P$ where all
Schubert basis elements take positive values has a
parameterization $(q_1,\dotsc, q_k) :
X_{P,>0}\overset\sim\To\R_{>0}^k$ given by the quantum
parameters. Secondly, this set lies in the smooth locus of $X_P$
and the inverse of the map giving the parameterization is
analytic. Thirdly, this set of Schubert positive points agrees
with the set of totally nonnegative matrices in $X_P$.
The Proposition~\ref{t:2} stated in the introduction is proved
immediately as a corollary.

In Section~\ref{s:Grassmannian} we make an excursion to recall
what these results look like explicitly in the Grassmannian case,
which is studied in detail in an earlier paper. We then use the
Grassmannian components of the Peterson variety to prove that the
top Schubert class $\si_{{w_0}^P}$ is generically nonvanishing as
function on $X_P$. Conjecturally, the same should hold for the
quantum Euler class, $\sum_{w\in W^P}\sigma_w\sigma_{PD(w)}$,
which would imply that $qH^*(G/P)$ is reduced.

The rest of the paper is devoted to the proof of the
Theorem~\ref{t:main}. In the next two sections, parts (1) and (2)
of the main theorem are proved.  The main ingredient for
constructing and parameterizing the Schubert positive points is
the positivity of the structure constants (Gromov--Witten
invariants). Computing the fiber in $X_P$ over a fixed positive
value of the quantum parameters $(q_1,\dotsc, q_k)$ is turned into
an eigenvalue problem for an irreducible nonnegative matrix, and
the unique Schubert positive solution we require is provided by a
Perron-Frobenius eigenvector. The smoothness property turns out to
be related to the positivity of the quantum Euler class, while
bianalyticity comes as a consequence of the one-dimensionality of
the Perron-Frobenius eigenspace.

The final part of Theorem~\ref{t:main}, that the notion of
positivity coming from Schubert bases agrees with total
positivity, is perhaps the most surprising. The problem is that,
except in the case of Grassmannian permutations, we know no useful
way to compute the Schubert classes as functions on the $X_P$. In
Section~\ref{s:SchubertClasses} we begin to simplify this problem
by proving another remarkable component of Peterson's theory.
Namely, consider the functions on $X_B$ given by the Schubert
classes of $qH^*(G/B)$. Then when extended as rational functions
to all of $X$, these restrict to give all the Schubert classes on
the smaller strata $X_P$. We prove this explicitly using quantum
Schubert polynomials and Fomin, Gelfand and Postnikov's quantum
straightening identity.

This last result enables us essentially to reduce the proof of the
final part of Theorem~\ref{t:main} to the full flag variety case.
The main problem there is to prove that an arbitrary Schubert
class takes positive values on the totally positive part. This is
done by topological arguments, using that the totally positive
part of $X_B$ is a semigroup.

\vskip.2cm {\it Acknowledgements.} Dale Peterson's beautiful
results presented by him in a series of lectures at MIT in 1997
were a major source of inspiration and are the foundation for much
of this paper. It is a pleasure to thank him here. Some of this
work was done during a very enjoyable and fruitful stay at the
Erwin Schr\"odinger Institute in Vienna, and sincere thanks go to
Peter Michor for his hospitality. Finally, I would like to thank
Bill Fulton for his kind invitation to Michigan, where some of the
final writing up could be done.

\section{Preliminaries}\label{s:Prelims}
Let $G=GL_n(\C)$, and $I=\{1,\dotsc, n-1\}$ an indexing set for
the simple roots. Denote  by $\Ad$  the adjoint representation of
$G$ on its Lie algebra $\mathfrak g$. We fix the Borel subgroups
$B^+$ of upper-triangular matrices and $B^-$ of lower-triangular
matrices in $G$. Their Lie algebras are denoted by $\mathfrak
b^+$ and $\mathfrak b^-$ respectively. We will also consider their
unipotent radicals $U^+$ and $U^-$ with their Lie algebras
$\mathfrak u^+$ and $\mathfrak u^-$ and the maximal torus
$T=B^+\cap B^-$. Let $X^*(T)$ be the character group of $T$ and
$X_*(T)$ the group of cocharacters with the usual perfect pairing
$<\, ,\, >: X^*(T)\x X_*(T)\to\Z$ between them. Let
$\Delta_+\subset X^*(T)$ be the set of positive roots
corresponding to $\mathfrak b^+$, and $\Delta_-$ the set of
negative roots. The fundamental weights and coweights are denoted
by $\omega_1,\dotsc, \omega_{n-1}\in X^*(T)$ and
$\omega^\vee_1,\dotsc,\omega^\vee_{n-1}\in X_*(T)$ respectively.
We define $\Pi\subset X^*(T)$ to be the set of positive simple
roots. The elements of $\Pi$ are denoted $\alpha_1,\dotsc,
\alpha_{n-1}$ where the $\alpha_m$-root space $\mathfrak
g_{\alpha_m}\in\mathfrak g$ is spanned by
\begin{equation*}
e_{\al_m}=\left(\delta_i^m\delta_j^{m+1}\right)_{i,j=1}^n=\begin{pmatrix}
 0& & & & &  \\
  &\ddots & & & &  \\
  & &0& 1& &  \\
  & && 0& & \\
  & & &  &\ddots &  \\
  & & &  & & 0
\end{pmatrix}.
\end{equation*}
Let $e :=\sum_{m=1}^{n-1} e_{\alpha_m}$. A special role will be
played by the principal nilpotent element $f\in \mathfrak u^-$
which is the transpose of $e$.

We identify the Weyl group $W$ of $G$, the symmetric group,  with
the group of all permutation matrices. $W$ is generated by the
usual simple reflections (adjacent transpositions) $s_1,\dotsc,
s_{n-1}$. The length function $\ell: W\to \N$ gives the length of
a reduced expression of $w\in W$ in the simple generators. There
is a unique longest element which is denoted $w_0$.

\subsection{Parabolics}\label{s:parabolics}
Let $P$ always denote a parabolic subgroup of $G$ containing $B^-$
and $\mathfrak p$ the Lie algebra of $P$. Let $I_P$ be the subset
of $I$ associated to $P$ consisting of all the $i\in I$ with
$s_i\in P$ and consider its complement $I^P:=I\setminus I_P$. We
will denote the elements of $I^P$ by $\{n_1,\dotsc, n_k\}$ with
\begin{equation*}
n_0:=0<n_1<n_2<\dotsc<n_k<n_{k+1}:=n.
\end{equation*}
Then the homogeneous space $G/P$ may be identified with the
partial flag variety (of quotients)
\begin{equation*}
G/P=\mathcal {F}_{n_1,\dotsc, n_k}(\mathbb C^n)=
\{\C^n\twoheadrightarrow
V_k\twoheadrightarrow\cdots\twoheadrightarrow V_1\to 0\ |\ \dim
V_j=n_j
 \}.
\end{equation*}
Next introduce $W_P=\left<s_i\ |\ i\in I_P\right>$, the parabolic
subgroup of $W$ corresponding to $P$. And $W^P\subset W$, the set
of minimal coset representatives for $W/W_P$. An element $w$ lies
in $W^P$ precisely if for all reduced expressions
$w=s_{i_1}\cdots s_{i_m}$ the last index $i_m$ always lies in
$I^P$. We write $w^P$ or $w_0^P$ for the longest element in $W^P$,
while the longest element in $W_P$ is denoted $w_P$. For example
$w^B_0=w_0$ and $w_B=1$. Finally $P$ gives rise to a decomposition
\begin{equation*}
\Delta_+= \Delta_{P,+}\sqcup \Delta_{+}^P,
\end{equation*}
where $\Delta_{P,+}=\{\al\in\Delta_+ \ | \ <\al,
\omega_i^\vee>=0
\text{ all $i\in I^P$} \} $ and $\Delta_+^P$ is its complement. So for
example $\Delta_{B,+}=\emptyset$ and $\Delta^B_+=\Delta_+$.

\section{The quantum cohomology ring of $G/P$}

\subsection{The usual cohomology of $\mathbf{G/P}$ and its Schubert
basis~}\label{s:ClassCoh} For our purposes it will suffice to
take homology or cohomology with complex coefficients, so always
$H^*(G/P)=H^*(G/P,\C)$. By the well-known result of C. Ehresmann,
the singular homology of the partial flag variety $G/P$ has a
basis indexed by the elements $w\in W^P$ made up of the
fundamental classes of the Schubert varieties,
\begin{equation*}
\Omega^P_w:=\overline{(B^+wP/P)}\subseteq G/P.
\end{equation*}
Here the bar stands for (Zariski) closure. Let $\si^P_{w}\in
H^*(G/P)$ be the Poincar\'e dual class to $[\Omega^P_w]$. Note
that $\Omega^P_w$ has complex codimension $\ell(w)$ in $G/P$ and
hence $\sigma^P_w$ lies in $H^{2\ell(w)}(G/P)$. The set
$\{\si^P_{w}\ |\ w\in W^P\}$ forms a basis of $H^*(G/P)$ called
the Schubert basis. The top degree cohomology of $G/P$ is spanned
by $\si^P_{w_0^P}$ and we have
the Poincar\'e duality pairing
\begin{equation*}
 H^*(G/P)\x H^*(G/P)\To \C ,\qquad(\si,\mu)\mapsto \left<\si\cdot \mu\right>
\end{equation*}
which may be interpreted as taking $(\si,\mu)$ to the coefficient
of $\si^P_{w^P_0}$ in the basis expansion of the product $\si\cdot
\mu$. For $w\in W^P$ let $PD(w)\in W^P$ be the minimal length
coset representative in $w_0wW_P$. Then this pairing is
characterized by
\begin{equation*}
\left <\si_w\cdot \si_v\right>=\delta_{w, PD(v)}.
\end{equation*}

\subsection{Definition of the quantum cohomology ring $\mathbf{qH^*(G/P)}$}
The (small) quantum cohomology ring $qH^*(G/P)$ may be defined by
enumerating curves into $G/P$ with certain properties. This
description is responsible for its positivity properties and is
the one we will give here. For more general background there are
already many books and survey articles on the subject of quantum
cohomology, see e.g.
\cite{CoxKatz:QCohBook,FuPa:QCoh,Manin:QCohBook,McDSal:QCohBook}
and references therein.

Let $I^P=\{n_1,\dotsc, n_k\}$. Then as a vector space the quantum
cohomology of the partial flag variety $G/P$ is given by
\begin{equation*}
qH^*(G/P)=\C[q^P_1,\dotsc, q^P_k]\otimes_{\overset{}\C} H^*(G/P),
\end{equation*}
where $q_1^P,\dotsc, q_k^P$ are called the quantum parameters.
Consider the Schubert classes as elements of $qH^*(G/P)$ by
identifying $\si^P_w$ with $1\otimes \si^P_w$. We will sometimes
drop the superscript $P$'s from the notation for the Schubert
classes and the quantum parameters when there is no possible
ambiguity.

Now $qH^*(G/P)$ is a free $\C[q^P_1,\dotsc, q^P_k]$-module with
basis given by the Schubert classes $\sigma^P_w$. It remains to
give the structure constants
$\left<\si_u^P,\si_v^P,\si_w^P\right>_{\mathbf d}$ in
\begin{equation*}
\sigma^P_v \sigma^P_w=\sum_{\begin{smallmatrix}u\in W^P\\ \mathbf
d\in \mathbb N^k
\end{smallmatrix} }\left<\si_u^P,\si_v^P,\si_w^P\right>_{\mathbf
d}\ \mathbf q^\mathbf d\si^P_{PD(u)}
\end{equation*}
to define the ring structure on $qH^*(G/P)$. Here $\mathbf
q^\mathbf d$ is multi-index notation for $\prod_{i=1}^k
q_i^{d_i}$.

Consider the set $\mathcal M_\mathbf d$ of holomorphic maps
$\phi:\mathbb {CP}^1\to G/P$, such that
\begin{equation*}
\phi_*\left (\left [\mathbb {CP}^1\right]\right)=\sum_{i=1}^k
d_i\left[\Omega^P_{s_{n_i}}\right ].
\end{equation*}
$\mathcal M_\mathbf d$ can be made into a quasi-projective
variety of dimension equal to $\dim(G/P)-\sum
d_i(n_{i+1}-n_{i-1})$. To define
$\left<\si_u^P,\si_v^P,\si_w^P\right>_{\mathbf d}$ first
translate the Schubert varieties $\Omega^P_u,\Omega^P_v$ and
$\Omega^P_w$ into general position, say to
$\widetilde{\Omega^P_u},\widetilde{\Omega^P_v}$ and
$\widetilde{\Omega^P_w}$.
Now consider the set $\mathcal M_\mathbf d(u,v,w)$ of all maps
$\phi\in \mathcal M_\mathbf d$ such that
\begin{equation*}
\phi(0)\in \widetilde{\Omega^P_u},\quad
\phi(1)\in\widetilde{\Omega^P_v},\quad \text{and}\quad
\phi(\infty)\in\widetilde{\Omega^P_w}.
\end{equation*}
Then  $\mathcal M_\mathbf d(u,v,w)$ is finite if $\dim(G/P)-\sum
d_i(n_{i+1}-n_{i-1})=\ell(u)+\ell(v)+\ell(w)$ and one may set
\begin{equation*}
\left<\si_u^P,\si_v^P,\si_w^P\right>_{\mathbf d}=\begin{cases}\#
\mathcal M_\mathbf d(u,v,w)& \text{if }\dim(G/P)-\sum
d_i(n_{i+1}-n_{i-1})=\\
&\qquad\qquad\qquad\qquad\qquad\quad=\ell(u)+\ell(v)+\ell(w),\\
0&\text{otherwise.}\end{cases}
\end{equation*}
These quantities are $3$-point, genus $0$ Gromov--Witten
invariants. By looking at $\mathbf d=(0,\dotsc, 0)$ one recovers
the classical structure constants obtained from intersecting
Schubert varieties in general position. Therefore this
multiplicative structure is a deformation the classical cup
product.
We note that the structure constants by their definition are {\it
nonnegative} integers.

The quantum cohomology analogue of the Poincar\'e duality pairing
may be defined as the symmetric $\C[q^P_1,\dotsc,
q^P_k]$-bilinear pairing
\begin{equation*} qH^*(G/P)\x qH^*(G/P)\To
\C[q^P_1,\dotsc, q^P_k], \qquad (\si,\mu)\mapsto \left<\si\cdot
\mu\right>_{\mathbf q}
\end{equation*}
which takes $(\si,\mu)$ to the coefficient of
$\si^P_{w_0^P}$ in the Schubert basis expansion of the product
$\si\cdot \mu$.
In terms of the Schubert basis the quantum Poincar\'e duality
pairing on $qH^*(G/P)$ is given by
\begin{equation*}
\left <\si^P_w\cdot\si^P_v\right>_{\mathbf q}=\delta_{w, PD(v)},
\end{equation*}
where $v,w\in W^P$, and $PD:W^P\to W^P$ is the involution defined
in Section \ref{s:ClassCoh} (see e.g.
\cite{Cio:QCohPFl}~Lemma~6.1).

\subsection{Borel's Presentation of the
cohomology ring $\mathbf{H^*(G/P)}$~}\label{s:Borel}

Let $G/P$ be realized as variety of flags of quotients as in
Section \ref{s:parabolics},
\begin{equation*}
G/P=\mathcal {F}_{n_1,\dotsc, n_k}(\mathbb C^n)=\left
\{\C^n=V_{k+1}\twoheadrightarrow
V_k\twoheadrightarrow\cdots\twoheadrightarrow V_1\to 0\ |\ \dim
V_j=n_j \right\}.
\end{equation*}
Then for $1\le j\le k+1$, the successive quotients $Q_j=\ker
(V_{j}\to V_{j-1})$ define rank $(n_{j}-n_{j-1})$ vector bundles
on $G/P$. Their Chern classes shall be denoted
\begin{equation*}
c_i(Q_j)=:\si^{(j)}_i=\si^{(j)}_{i,P}.
\end{equation*}
By the splitting principle it is natural to introduce independent
variables $x_1,\dotsc, x_n$ such that $x_{n_{j-1}+1},\dotsc,
x_{n_{j}}$ are the Chern roots of $Q_j$. So
\begin{equation}\label{e:sigmas}
\si^{(j)}_{i}=e_i(x_{n_{j-1}+1},x_{n_{j-1}+2},\dotsc ,x_{n_j})
\end{equation}
the $i$-th elementary symmetric polynomial in the variables
$\{x_{n_i+1},\dotsc, x_{n_{i+1}}\}$. Let $W_P$ act on the
polynomial ring $\C[x_1,\dotsc, x_n]$ in the natural way by
permuting the variables. Then the ring of invariants is
\begin{equation*}
\C[x_1,\dotsc,x_n]^{W_P}=\C\left
[\si^{(1)}_1,\dotsc,\si^{(1)}_{n_1},\si^{(2)}_1,
\dotsc\dotsc,\si^{(k+1)}_1,\dotsc,\si^{(k+1)}_{n-n_k}\right ].
\end{equation*}
A. Borel \cite{Borel:CohG/P} showed that the Chern classes
$\si^{(j)}_{i}$ generate $H^*(G/P)$ and the only relations
between these generators come from the triviality of the bundle
$Q_1\oplus Q_2\oplus\cdots \oplus Q_k$ (which may be trivialized
using a Hermitian inner product on $\C^n$).

In other words, if $J$ denotes the ideal in
$\C[x_1,\dotsc,x_n]^{W_P}$ generated by the elementary symmetric
polynomials $e_i(x_1,\dotsc, x_n)$, then
\begin{equation}\label{e:BorelIso}
H^*(G/P)\cong\C[x_1,\dotsc,x_n]^{W_P}/J.
\end{equation}

\subsection{Schubert polynomials and elementary monomials for $\mathbf
{H^*(G/P)}$}\label{s:Schubs}

For $1\le i\le k+1$, define
\begin{equation}\label{e:ElSymPol}
e_i^{(j)}=e_{i,P}^{(j)}:=e_i(x_1,\dotsc, x_{n_j}),
\end{equation}
the $i$-th elementary symmetric polynomial in $n_j$ variables.
Then the $e^{(k+1)}_i$'s are the generators of the ideal $J$. But
for $1\le j\le k$, the element $e_i^{(j)}$ corresponds under
\eqref{e:BorelIso} to a nonzero element of $ H^*(G/P)$, namely the
special Schubert class $\sigma^P_{s_{n_j-i+1}\cdots s_{n_j}}$.

The polynomial defined by
\begin{equation}\label{e:topclass}
c_{w_0^P}:=\left(e^{(1)}_{n_1}\right )^{n_1}\cdots
\left(e^{(k)}_{n_k}\right)^{n_k-n_{k-1}}
\end{equation}
represents the top class $\sigma_{w_0^P}$.

These are examples of the {\it Schubert polynomials} of Lascoux and
Sch\"utzenberger, \cite{LascSch:SchubPol}. The
Schubert polynomials $\{c_w \ |\ w\in W^P \}\subset
\C[x_1,\dotsc,x_n]^{W_P}$ are, loosely speaking, obtained from
the top one by divided difference operators, see
\cite{LascSch:SchubPol}
or \cite{Macd:SchubPol} for details. If $P=B$ then Schubert
polynomials $c_w$ are obtained for all $w\in W$, and
the ones from above corresponding to $G/P$ are just the subset consisting of
all those for which $w\in W^P$. The key property of a Schubert
polynomial $c_w$ is of
course that it is a representative for the corresponding Schubert
class $\sigma_w$.

A different description, following \cite{FoGePo:QSchub}, of the
Schubert polynomials $c_w$ for $w\in W^P$ says precisely where
these representatives must lie. They are those representatives of
the Schubert classes which may be written as linear combinations
of certain ``elementary monomials'' in $\C[x_1,\dotsc,
x_n]^{W_P}$.

Explicitly, let $\mathcal L_P$ be the set of sequences
$\Lambda=(\lambda^{(1)},\dotsc,\lambda^{(k)})$ of partitions, such
that $\lambda^{(j)}$ has at most $(n_j-n_{j-1})$ parts and
$\lambda^{(j)}_1\le n_j$. To any $\Lambda\in\mathcal L_P$
associate a polynomial,
\begin{equation*}
e_{\Lambda}=\left (e^{(1)}_{\lambda^{(1)}_1}\cdots
e^{(1)}_{\lambda^{(1)}_{n_1}}\right)\cdots
\left(e^{(k)}_{\lambda^{(k)}_1}\cdots
e^{(k)}_{\lambda^{(k)}_{n_k-n_{k-1}}}\right).
\end{equation*}
Let us call these polynomials {\it $P$-standard monomials}. These
$e_\Lambda$ are linearly independent and span a complementary
subspace to the ideal $J$. So
\begin{equation*}
\C[x_1,\dotsc, x_n]^{W_P}=J\ \oplus \left <e_\Lambda\right
>_{\Lambda\in \mathcal L_P}.
\end{equation*}
Then the Schubert polynomial $c_{w}$ is the (unique) representative in $\left
<e_\Lambda\right >_{\Lambda\in \mathcal L_P}$ for the Schubert class $\si^P_{w}$.

\subsection{Astashkevich, Sadov and Kim's
presentation of $\mathbf {qH^*(G/P)}$} The presentation for the
quantum cohomology ring $qH^*(G/P)$ analogous to Borel's
presentation of $H^*(G/P)$ was first discovered by Astashkevich
and Sadov \cite{AstSa:QCohPFl} and Kim \cite{Kim:QCohPFl}. A
complete proof may be found in Ciocan-Fontanine
\cite{Cio:QCohPFl}. In special cases these presentations were
known earlier, e.g. for Grassmannians
\cite{Be:QSchubCalc,SiTi:QCoh},
and in the full flag variety case \cite{GiKi:FlTod,Cio:QCohFl}.

The generators of $qH^*(G/P)$ will be the generators of the
usual cohomology ring $\si^{(j)}_{i}$ (embedded as
$1\otimes\si^{(j)}_{i}$) along with the quantum parameters
$q^{P}_1,\dotsc, q^{P}_k$. Here $i,j$ runs through $1\le j\le k+1$
and $1\le i\le n_j-n_{j-1}$.
Let us for now treat the $\si^{(j)}_i$ and $q_j$ as independent
variables generating a polynomial ring
$\C[\si^{(1)}_1,\dotsc,\si^{(k+1)}_{n-n_k}\ , \ q_1,\dotsc, q_k]$.

\begin{defn}[$(\mathbf q,P)$-elementary symmetric
polynomials]\label{d:E} Let $i\in\Z$ and
$l\in\{-1,0,\dotsc,k+1\}$. Define elements
$E^{(l)}_{i,P}=E^{(l)}_{i}\in\C[\si^{(1)}_1,\dotsc,\si^{(k+1)}_{n-n_k}\,
,\,  q_1,\dotsc,q_k]$ recursively as follows. The initial values
are
\begin{equation*}
E^{(-1)}_i= E^{(0)}_i=0\ \text{ for all $i$},\  \ \text{and }\
E^{(l)}_{i}=0 \ \text{ unless } 0\le i\le n_l ,
\end{equation*}
and we set $\si^{(l)}_i=0$ if $i>n_l-n_{l-1}$ and
$\sigma^{(l)}_0=1$. For $1\le l\le k+1$ and $0\le i\le n_l$ the
polynomial $E^{(l)}_i$ is defined by
\begin{multline*}
E^{(l)}_i\hskip -.1cm=\hskip -.1cm\left (E^{(l-1)}_{i} +
\si^{(l)}_1 E^{(l-1)}_{i-1}+\cdots+ \si^{(l)}_{i-1} E^{(l-1)}_{1}+
\si^{(l)}_{i} \right ) +\\ (-1)^{n_{i+1}-n_{i}+1}q_{l-1}
E^{(l-2)}_{i-n_{l}+n_{l-2}}.
\end{multline*}
If the $q_l$ are set to $0$ and the $\si^{(l)}_i$ are as in
\eqref{e:sigmas}, then this recursion defines the elementary
symmetric polynomials $e^{(l)}_i$.
\end{defn}

\begin{rem} This is a basic recursive definition of the
quantum elementary symmetric polynomials. See \cite{Cio:QCohPFl}
for a host of other descriptions. And here is also one other curious one
to add to this list.

Order the variables $\sigma_i^{(j)}$ lexicographically, so that
$\sigma_i^{(j)}<\sigma_k^{(j')}$ whenever $j<j'$ or $j=j'$ and
$i<k$. Now suppose just for the remainder of this remark that the
variables $\sigma_i^{(j)}$ are not necessarily commuting. More
precisely, let $\sigma_i^{(j)}$ and $\sigma_k^{(j')}$ commute unless
$|j-j'|=1$ and both $i$ and $k$ are maximal. In that case impose
the Heisenberg relation

\begin{equation*}
\sigma_{n_j-n_{j-1}}^{(j)}\sigma_{n_{j+1}-n_{j}}^{(j+1)}=
\sigma_{n_{j+1}-n_{j}}^{(j+1)}
\sigma_{n_j-n_{j-1}}^{(j)} + q_j.
\end{equation*}
The $q_i$ commute with everything. Now add a central variable $x$
and define polynomials $p_j(x)=x^{n_j-n_{j-1}}+
\si^{(j)}_{1}x^{n_j-n_{j-1}-1}+\dotsc +
\si^{(j)}_{n_{j}-n_{j-1}}$. Then expand the product
\begin{equation*}
p_m(x)\cdot p_{m-1}(x)\cdot\dotsc\cdot p_1(x)
\end{equation*}
and write the resulting coefficients in terms of increasing
monomials in the $\sigma_i^{(j)}$ (monomials with factors ordered
in increasing fashion) by using the commutation relations. Then
for $d\le n_m$ the coefficient of $x^{n_m-d}$ gives the $(\mathbf
q,P)$-elementary symmetric polynomial $E^{(m)}_d$.

For example in type $A_2$ for the full flag variety case, the
polynomials $E^{(3)}_1$, $E^{(3)}_2$ and $E^{(3)}_3$ turn up as
coefficients in
\begin{multline*}
(x+x_3)(x+x_2)(x+x_1)=x^3+(x_1+x_2+x_3)x^2+\\(x_1x_2+x_2
x_3+x_1x_3+q_1+q_2)x + (x_1 x_2 x_3 + x_1 q_2+ x_3 q_1).
\end{multline*}
\end{rem}

\begin{thm}\cite{Cio:QCohPFl,AstSa:QCohPFl,Kim:QCohPFl}\label{t:Cio}
The assignment $\si^{(i)}_j\mapsto 1\otimes\si^{(i)}_j$ and
$q_i\mapsto q_i\otimes 1$ gives rise to an isomorphism
\begin{equation}\label{e:presentation}
\C[\,\si^{(1)}_{1},\dotsc,\si^{(k)}_{n-n_k}\, ,\,
q_{1},\dotsc,q_{k}\,]/J \overset\sim\longrightarrow qH^*(G/P),
\end{equation}
where $J$ is the ideal $(E^{(k+1)}_{1},\dotsc,E^{(k+1)}_{n})$.
This isomorphism takes the element $E^{(l)}_i$ for $1\le l\le k$
to the special Schubert class $\si^P_{s_{n_l-i+1}\cdots s_{n_l}}$.
\end{thm}
An immediate question raised by this theorem is how to describe
Schubert classes in the picture on left hand side of this
isomorphism. This is answered by a quantum analogue of the
Schubert polynomials.

\subsection{Quantum Schubert Polynomials}\label{s:Schub} In the
case of $G/B$ a full theory of quantum Schubert polynomials was
given by Fomin, Gelfand and Postnikov \cite{FoGePo:QSchub}. This
was later generalized to partial flag varieties by
Ciocan-Fontanine \cite{Cio:QCohPFl}. Note that the quantum
Schubert polynomials for partial flag varieties are {\it not}
special cases of the full flag variety ones, due to lack of
functoriality of quantum cohomology (but see Proposition
\ref{p:SchubRestr}).

There is also a different construction of (double) quantum
Schubert polynomials due to Kirillov and Maeno
\cite{KirANMa:QSchub} which has been shown to give the same
answer. We give the definitions following \cite{FoGePo:QSchub} and
\cite{Cio:QCohPFl} below.

\begin{defn}[$(\mathbf q,P)$-standard monomials
] As in Section~\ref{s:Schubs} let $\mathcal L_P$ be the set of
sequences $\Lambda=(\lambda^{(1)},\dotsc,\lambda^{(k)})$ of
partitions, such that $\lambda^{(j)}$ has at most $(n_j-n_{j-1})$
parts and $\lambda^{(j)}_1\le n_j$. To each $\Lambda\in\mathcal
L_P$ associate an element
\begin{equation*}
E_{\Lambda}=\left (E^{(1)}_{\lambda^{(1)}_1}\cdots
E^{(1)}_{\lambda^{(1)}_{n_1}}\right)\cdots
\left(E^{(k)}_{\lambda^{(k)}_1}\cdots
E^{(k)}_{\lambda^{(k)}_{n_k-n_{k-1}}}\right)\in
\C[\si^{(1)}_1,\dotsc,\si^{(k)}_{n-n_k},q_1,\dotsc, q_k].
\end{equation*}
These elements are called the {\it $(\mathbf q,P)$-standard
monomials}.
\end{defn}

\begin{rem}[Example] The $(\mathbf
q,B)$-standard polynomials were introduced in
\cite{FoGePo:QSchub}. They are the monomials of the form
\begin{equation*}
E^{(1)}_{B,j_1}E^{(2)}_{B,j_2}\cdots E^{(n-1)}_{B,j_{n-1}}
\end{equation*}
where $0\le j_l\le l$ for all $l=1,\dotsc,n-1$.
\end{rem}

Let $V$ denote the $\C[q_1,\dotsc, q_k]$-module
spanned by the $(\mathbf q,P)$-standard monomials,
\begin{equation*}
V=\C[q_1,\dotsc, q_k]\otimes_\C \left<E_\Lambda\right
>_{\Lambda\in\mathcal L_P}.
\end{equation*}
Then
\begin{equation*}
\C[\si^{(1)}_1,\dotsc,\si^{(k)}_{n-n_k},q_1,\dotsc, q_k]=J\
\oplus\  V.
\end{equation*}
\begin{defn}[$(\mathbf q,P)$-Schubert polynomials]
The {\it quantum Schubert polynomial}
$C^P_{w}\in\C[\si^{(1)}_1,\dotsc,\si^{(k+1)}_{n-n_k},q_1,\dotsc,
q_k]$ is defined to be the unique element of $V$ whose coset
modulo $J$ maps to the Schubert class $\si^P_{w}$ under the
isomorphism
\begin{equation*}
\C[\si^{(1)}_1,\dotsc,\si^{(k+1)}_{n-n_k},q_1,\dotsc,
q_k]/J\overset\sim\To qH^*(G/P).
\end{equation*}
\end{defn}

\begin{rem}[Example] From \eqref{e:topclass} it follows immediately that the
$(\mathbf q,P)$-Schubert polynomial representing the top class
$\si_{w_0^P}\in qH^*(G/P)$ is given by
\begin{equation}\label{e:qtopclass}
C^P_{w_0^P}=\left(E^{(1)}_{P,n_1}\right )^{n_1}\cdots
\left(E^{(k)}_{P,n_k}\right)^{n_k-n_{k-1}}.
\end{equation}
\end{rem}

\subsection{Grassmannian permutations}\label{s:Kirillov}

A Grassmannian permutation of descent $m$ is an element $w\in
W^{P_d}$ for the maximal parabolic $P_d$ with $I^{P_d}=\{d\}$. As
permutations on $\{1,\dotsc, n\}$ these may be characterized by
\begin{equation*}
w\in W^{P_d}\ \iff\ w(1)<\cdots<w(d)\ \text{ and }\ w(d+1)<\cdots
<w(n).
\end{equation*}
They are in bijective correspondence with shapes (partitions)
$\lambda=(\lambda_1,\dotsc, \lambda_d)$ such that $n-d\ge
\lambda_1\ge \dotsc \ge \lambda_d\ge 0$, via $\lambda_i=w(i)-i+1$.

Let $w_{\lambda,d}$ denote the Grassmannian permutation of
descent $d$ and shape $\lambda$. There is a closed formula for the
quantum Schubert polynomials $C^B_{w_{\lambda,d}}$ given by A. N.
Kirillov in \cite{KirAN:QSchurFcts} which we will derive here from
Fomin, Gelfand and Postnikov's definition.

The classical Schubert polynomial for $w_{\lambda,d}$ is just the
Schur polynomial $c_{w_{\lambda,d}}=s_\lambda(x_1,\dotsc, x_d)$,
see e.g. \cite{Macd:SchubPol}. Therefore by the Jacobi-Trudi
identity
\begin{equation*}
c_{w_{\lambda,d}}=\det
\begin{pmatrix}
 e_{\lambda'_1}^{(d)} & e_{\lambda'_1+1}^{(d)} &\cdots &   e_{\lambda'_1+c-1}^{(d)} \\
 e_{\lambda'_2-1}^{(d)} &e_{\lambda'_2}^{(d)} & \cdots  & e_{\lambda'_2+c-2}^{(d)}  \\
               &               &  \ddots &                  \\
 e_{\lambda'_c+c-1}^{(d)}   &    \cdots           & \cdots &e_{\lambda'_{c}}^{(d)}
\end{pmatrix},
\end{equation*}
where $\lambda'$ is the conjugate partition to $\lambda$ and
$c=n-d$ (see \cite{Macd:SFBook}). Repeatedly applying the identity
\begin{equation*}
e^{(m)}_j=e^{(m+1)}_j - x_{m+1} e^{(m)}_{j-1}
\end{equation*}
of elementary symmetric polynomials to one column in the
determinant at a time, one obtains
\begin{equation*}
c_{w_{\lambda,d}}=\det
\begin{pmatrix}
 e_{\lambda'_1}^{(d)} & e_{\lambda'_1+1}^{(d+1)} &\cdots &   e_{\lambda'_1+c-1}^{(n-1)} \\
 e_{\lambda'_2-1}^{(d)} &e_{\lambda'_2}^{(d+1)} & \cdots  & e_{\lambda'_2+c-2}^{(n-1)}  \\
               &               &  \ddots &                  \\
 e_{\lambda'_c+c-1}^{(d)}   &    \cdots           & \cdots &e_{\lambda'_{c}}^{(n-1)}
\end{pmatrix}.
\end{equation*}
Expanding out this determinant gives an expression for
$c_{w_{\lambda,d}}$ as linear combination of  $B$-standard
monomials. Therefore the quantization is simply given by
\begin{equation*}
C_{w_{\lambda,d}}=\det
\begin{pmatrix}
 E_{\lambda'_1}^{(d)} & E_{\lambda'_1+1}^{(d+1)} &\cdots &   E_{\lambda'_1+c-1}^{(n-1)} \\
 E_{\lambda'_2-1}^{(d)} &E_{\lambda'_2}^{(d+1)} & \cdots  & E_{\lambda'_2+c-2}^{(n-1)}  \\
               &               &  \ddots &                  \\
 E_{\lambda'_c+c-1}^{(d)}   &    \cdots           & \cdots &E_{\lambda'_{c}}^{(n-1)}
\end{pmatrix}.
\end{equation*}

\subsection{Quantum Chevalley Formula}\label{s:QChev}
The Pieri formula for $H^*(G/P)$ of Lascoux and Sch\"utzenberger
was generalized to the $qH^*(G/P)$ setting by Ciocan-Fontanine in
\cite{Cio:QCohPFl}. We will only need the following simpler case.

\begin{thm}
\label{t:qChev} For $h\le l\in \{1,\dotsc, k\}$ set
$\tau_{h,l}=s_{n_h}\cdot\dotsc\cdot s_{n_{l+1}-1}s_{n_l-1}\cdot
\dotsc\cdot s_{n_{h-1}+1}$ and $\mathbf q_{h,l}=q_h\cdot q_{h+1}
\dotsc \cdot q_l$.  Let $n_j\in I^P$ and $w\in W^P$. Then
\begin{equation*}
\si^P_{s_{n_j}}\si^P_{w}=
\sum_{
\begin{smallmatrix}\alpha\in\Delta^+\\
ws_\alpha\in W^P\\
\ell(ws_\alpha)=\ell(w)+1
\end{smallmatrix}}
<\alpha,\omega_{n_j}^\vee>\si^P_{w s_\alpha} + \sum_{
\begin{smallmatrix}
h,l\in\{1,\dotsc, k\}\\
 1\le h\le j\le l \le k\\
\ell(w\tau_{h,l})=\ell(w)-\ell(\tau_{h,l})
\end{smallmatrix}}
\mathbf q_{h,l}\si^P_{w\tau_{h,l}}.
\end{equation*}
\end{thm}
\noindent This is a reformulation of a special case of Theorem 3.1
in \cite{Cio:QCohPFl}.

\section{Quantum cohomology rings as coordinate rings}\label{s:Peterson}

\subsection{ASK-matrices}\label{s:ASKmatrices} We introduce with
some minor changes an $n\x n$ matrix $A^{[k+1]}$ with entries in
$\C[\si^{(1)}_1,\dotsc, \si^{(k)}_{n-n_k}\, ,\, q_1,\dotsc, q_k]$
introduced by Astashkevich, Sadov and Kim in \cite{AstSa:QCohPFl}
and \cite{Kim:QCohPFl}. Setting $n_{k+1}=n$ and $n_0=0$ define
first $(n_{j}-n_{j-1})\x (n_{j}-n_{j-1})$ matrices $D^{(j)}$ by
\begin{equation*}
D^{(1)}=\begin{pmatrix} -\si^{(1)}_1 & -\si^{(1)}_2& \dotsc &
 -\si^{(1)}_{n_1}\\
 0 &\cdots&\cdots&0 \\
 \vdots&&&\vdots\\
  \vdots&&&\vdots\\
 0&\cdots &\cdots&0
\end{pmatrix}\ \text {and}\ \
D^{(j)}=\begin{pmatrix}
0&\cdots&0&-\si^{(j)}_{n_j-n_{j-1}}\\
 \vdots&&\vdots&\vdots\\
 \vdots&&\vdots&-\si^{(j)}_2 \\
 0 &  \dotsc&0&
 -\si^{(j)}_{1}
\end{pmatrix}
\end{equation*}
for $2\le j\le k+1$. And let $D^{[l]}$ be the $n_l\x n_l$ block
matrix made up of diagonal blocks $D^{(1)},\dotsc, D^{(l)}$.
Furthermore define $n_l\x n_l$ matrices
\begin{equation*}
f^{[l]}=\begin{pmatrix}
 0& & &  \\
 1 &\ddots & &   \\
  &\ddots &\ddots&  \\
  &  &1 &0  \\
\end{pmatrix}\ \text{and}\ \ Q^{[l]}:=\left
((-1)^{n_{m+1}-n_m}q_m\delta_i^{n_{m-1}+1}\delta_j^{n_{m+1}}\right)_{i,j=1}^{n_l}.
\end{equation*}
Then set
\begin{equation*}
A^{[l]}:=f^{[l]}+D^{[l]}+Q^{[l]}.
\end{equation*}
The coefficients of the characteristic polynomials of the
$A^{[l]}$ satisfy precisely the same recursion as the $(P,\mathbf
q)$-standard symmetric polynomials $E^{(l)}_{P,i}$. Explicitly,
we have
\begin{equation*}
\det(\lambda \Id - A^{[l]})=\lambda^{n_l}+
E^{(l)}_{P,1}\lambda^{n_l-1}+\dotsc + E^{(l)}_{P,n_l}.
\end{equation*}
In particular the relations $E^{(k+1)}_1=\dotsc=E^{(k+1)}_n=0$ of
the quantum cohomology ring are equivalent to
\begin{equation}\label{e:nilpASK}
\det(\la\Id-A^{[k+1]})=\la^{n}.
\end{equation}
Let us call the matrices in $\mathfrak{gl}_n$ of the same form as
$A^{[k+1]}$ (with the same pattern of $0$ and $1$ entries) {\it
ASK-matrices}. They form an affine subspace $\mathcal A_P$ in
$\mathfrak {gl}_n$. Let $\mathcal N_P$ be the (non-reduced)
intersection,
\begin{equation*}
\mathcal N_P=\mathcal A_P\cap \mathcal N,
\end{equation*}
of $\mathcal A_P$ with the nilpotent cone $\mathcal N$ in
$\mathfrak {gl}_n$. Its coordinate ring is denoted $\mathcal
O(\mathcal N_P)$.

Then \eqref{e:nilpASK} implies that the map $\mathcal O(\mathcal
A_P) \To \C[\,\si^{(1)}_{1},\dotsc,\si^{(k)}_{n-n_k}\, ,\,
q_{1},\dotsc,q_{k}\,]$ defined by $A^{[k+1]}$ induces an
isomorphism
\begin{equation}\label{e:ON_Ppres}
\mathcal O(\mathcal N_P) \overset \sim\To
\C[\,\si^{(1)}_{1},\dotsc,\si^{(k)}_{n-n_k}\, ,\,
q_{1},\dotsc,q_{k}\,]/J.
\end{equation}
The statement of Theorem~\ref{t:Cio} may therefore be interpreted
as
\begin{equation}\label{e:ON_P}
\mathcal  O(\mathcal N_P)\overset\sim\To qH^*(G/P).
\end{equation}

\subsection{Peterson's Theorem}\label{s:Pet}
All the affine varieties $Spec(qH^*(G/P))$ turn out to be most
naturally viewed as embedded in the flag variety (or in general
in the Langlands dual flag variety) where they patch together as
strata of one remarkable projective variety called the Peterson
variety. This is the content of Dale Peterson's theorem which we
will deduce here explicitly for the type $A$ case.

\begin{defn}
Let $\mathfrak b^{++}:=\sum_{\alpha\in\Delta_+\setminus\Pi} \,
\mathfrak g_{\alpha}$ and $\pi^{++}:\mathfrak g\to\mathfrak
b^{++}$ is the projection along weight spaces. Let $f\in
\mathfrak{gl}_n$ be the principal nilpotent $f^{[k+1]}$ from
above. Then the equations
\begin{equation*}
\pi^{++}(\Ad(g\inv)\cdot f)=0
\end{equation*}
define a closed subvariety of $G$ invariant under right
multiplication by $B^-$. Thus they define a closed subvariety of
$G/B^-$. This subvariety $\Y$ is the {\it Peterson variety} for
type $A$. Loosely, $\Y$ can be described by
\begin{equation*}
\Y=\left\{gB^-\in G/B^-\ \left |\ \Ad(g\inv)\cdot f\in\mathfrak
b^- \oplus \sum_{i\in I} \C e_{\al_i}\right.\right\}.
\end{equation*}
\end{defn}

Let $v_1,\dotsc, v_n$ be the standard basis of $V=\C^n$. Then
$\{v_{i_1}\wedge\cdots\wedge v_{i_j}\ |\
0<i_1<i_2<\cdots<i_j<n\}$ is the standard basis of the
fundamental representation $V^{\omega_j}=\bigwedge^{j} V$. Let
$(\ \  |\ \ )$ denote the inner product on $V^{\omega_j}$ such
that this basis is orthonormal. We also refer to representations
by their lowest weight, so $V^{\omega_{n-m}}=:V_{-\omega_m}$.
Denote the lowest weight vector by
\begin{equation*}
v_{-\omega_m}=v_{m+1}\wedge\cdots\wedge v_n.
\end{equation*}
Define rational functions $G^{m}_i=G_{s_{m-i+1}\cdots s_{m}}$ on
$G/B^-$ in terms of matrix coefficients of the fundamental
representations by
\begin{equation}\label{e:Gs}
G^{m}_i(gB^-)=G_{s_{m-i+1}\cdots s_{m}}(gB^-):= \frac{(g\cdot
v_{-\omega_m}\ |\ s_{m-i+1}\cdots s_{m}\cdot v_{-\omega_m}) }
{(g\cdot v_{-\omega_m}\ |\ v_{-\omega_m})}.
\end{equation}
Note that $s_{m-i+1}\cdots s_{m}\cdot v_{-\omega_m}=
v_{m-i+1}\wedge v_{m+2}\wedge \cdots\wedge v_{n}$ and
$G^m_i(gB^-)$ may be written down simply as a quotient of two
$(n-m)\x (n-m)$-minors of $g$.

\begin{thm}[D. Peterson]\label{t:Pet}
\begin{enumerate}
\item For any parabolic subgroup $W_P\subset W$ with longest element
$w_P$ define $\Y_P$ as (non-reduced) intersection by
\begin{equation*}
\Y_P:= \Y\cap \left(B^+w_P B^-/B^-\right ).
\end{equation*}
Then on points these give a decomposition
\begin{equation*}
\Y(\C)=\bigsqcup_{P} \Y_P(\C).
\end{equation*}
\item
For each parabolic $P$ there is a unique isomorphism
\begin{equation}\label{e:PetIso1}
\mathcal O(\Y_P)\overset\sim\To qH^*(G/P) ,
\end{equation}
which ends $G_{s_{n_j-i+1}\cdots s_{n_j}}$ to
$\si^P_{s_{n_j-i+1}\cdots s_{n_j}}$.
\end{enumerate}
\end{thm}

\begin{rem}\label{r:coords} If $P$ is the parabolic subgroup,
then $G^m_j$ is a well-defined (regular) function on the Bruhat cell
$B^+w_P B^-/B^-$
precisely if $m\in I^P=\{n_1,\dotsc, n_k\}$. In fact we have
\begin{align*}
&B^+w_PB^-/B^- \overset\sim\To \quad\C^{\left(\sum_{i=1}^k n_i\right)}\\
&gB^- \mapsto \ (G^{n_1}_1(gB^-),\dotsc,
G^{n_1}_{n_1}(gB^-),G^{n_2}_1(gB^-),\dotsc,G^{n_k}_{n_k}(gB^-)),
\end{align*}
or in other words,
\begin{equation*}
\mathcal O(B^+w_PB^-/B^-)=\C[G^{n_1}_1,\dotsc,G^{n_k}_{n_k}].
\end{equation*}
Let $\mathcal J_P\subset \C[G^{n_1}_1,\dotsc,G^{n_k}_{n_k}]$
denote the ideal defining $\Y_P$.
\end{rem}

\begin{proof}\label{p:Pet}
(1) is proved in \cite{Pet:QCoh}. See also Lemma 2.3 in
\cite{Rie:QCohGr}.
We will deduce (2) very explicitly from the ASK presentation.
Begin by defining a particular section $u:B^+w_PB^-/B^-\to U^+$
of the map $x\mapsto x w_P B^-$ in the other direction. For
$l=0,\dotsc, k$ we have $n\x (n_{l+1}-n_{l})$ matrices $U^{(l)}$
defined by
\begin{equation*}
 U^{(0)}=
\begin{pmatrix}
1 & G^{n_1}_1 & G^{n_1}_2 & \cdots & G^{n_1}_{n_1-1}\\
0 & 1         & G^{n_1}_1 & \ddots &  \vdots        \\
\vdots &      &\ddots     &        & G^{n_1}_2      \\
  &           &           & \ddots &  G^{n_1}_1     \\
  &           &           &        &  1             \\
  &           &           &        &  0             \\
  &           &           &        &  \vdots        \\
  &           &           &        &                \\
  &           &           &        &                \\
\vdots &      &           &        &  \vdots        \\
0 & \cdots    &           & \cdots &  0
\end{pmatrix},
\quad U^{(l)}=
\begin{pmatrix}
 G^{n_l}_{n_l} &            &           &                \\
 \vdots        &  \ddots    &           &                \\
 \vdots        &            &\ddots     &                \\
 G^{n_l}_2     &            &           &G^{n_l}_{n_l}   \\
 G^{n_l}_{1}   & \ddots     &           &  \vdots        \\
 1             & \ddots     & \ddots    &  \vdots        \\
 0             &  \ddots    & G^{n_l}_1 & G^{n_l}_2      \\
 \vdots        &            & 1         &G^{n_l}_{1}     \\
               &            &           & 1              \\
 \vdots        &            &           &                \\
 0             &    \cdots  &\cdots     & 0
\end{pmatrix}
\end{equation*}
where $1\le l\le k$. Then the matrix
\begin{equation*}\label{e:u}
u=\left(\left.\left.\left.
\begin{matrix}
       &\\
          &  \\
       &\\
U^{(0)}& \\
       &\\
          &  \\
       &\\
\end{matrix}
\right |
\begin{matrix}
 &         &  \\
 &         &  \\
 &          & \\
 & U^{(1)} &  \\
  &         & \\
 &         &  \\
   &         &\\
\end{matrix}
\right |
\begin{matrix}
   &        & \\
  &         &  \\
 &        & \\
  & \cdots &  \\
 &         &  \\
  &        & \\
    &       & \\
\end{matrix}
\right |
\begin{matrix}
   &\\
   &\\
U^{(k)} \\
    &\\
   &\\
\end{matrix}
\right).
\end{equation*}
defines a map $u:B^+ w_P B^-/B^-\to U^+$. It follows using
Remark~\ref{r:coords} that this map is indeed a section. That is,
\begin{equation*}
gB^-=u(gB^-)w_P B^-, \quad \text{for\ \ $gB^-\in B^+w_P B^-/B^-$.}
\end{equation*}

Consider the matrix
\begin{equation*}
\tilde A=u\inv f u \quad \in \quad
\mathfrak{gl}_n(\C[G^{n_1}_1,\dotsc,G^{n_k}_{n_k}]).
\end{equation*}
A direct computation shows that modulo the ideal $\mathcal J_P$
defining $\Y_P$ the matrix $\tilde A$ is an ASK matrix (i.e. it
is an ASK matrix over
$\C[G^{n_1}_1,\dotsc,G^{n_k}_{n_k}]/\mathcal J_P$). Also it is
clear that the characteristic polynomial of $\tilde A$ satisfies
$\det(\lambda\Id - \tilde A)=\lambda^n$, since $\tilde A$ is
conjugate to $f$. Therefore the morphism $B^+ w_P B^-/B^-\To
\mathfrak{gl_n}$ defined by $\tilde A$ restricts to a morphism
\begin{equation}\label{e:Atilde}
\tilde A|_{\Y_P}: \Y_P\to \mathcal N_P
\end{equation}
from $\Y_P$ to the variety of nilpotent ASK-matrices.

For the inverse define a map $\psi:
\C[G^{n_1}_1,\dotsc,G^{n_k}_{n_k}]\to
\C[\si^{(1)}_1,\dotsc,\si^{(k+1)}_{n-n_k}\, ,\, q_1,\dotsc, q_k]$
by $\psi(G^{n_j}_i)= E^{(j)}_i$.
Applying $\psi$ to the entries of $u$ we obtain a matrix $u_E$
with entries in $\C[\si^{(1)}_1,\dotsc,\si^{(k+1)}_{n-n_k}\, ,\,
q_1,\dotsc, q_k]$. Then the recursive definition of the
$E^{(j)}_i$ translates into the identity
\begin{equation}\label{e:matrixid}
u_E\inv \, f\, u_E=A^{[k+1]}+ M
\end{equation}
of matrices over
$\C\left[\si^{(1)}_{1},\dotsc,\si^{(k+1)}_{n-n_k}\, ,\,
q_{1},\dotsc,q_{k}\right]$, where $A^{[k+1]}$ is the matrix
defined in Section \ref{s:ASKmatrices}, and $M$ the $n\x n$ matrix
given by
\begin{equation*}M=
\begin{pmatrix}
0&\cdots &\cdots & 0&E^{(k+1)}_n\\
\vdots& & & \vdots& E^{(k+1)}_{n-1}\\
\vdots& & &\vdots &     \vdots\\
\vdots& & & \vdots& E^{(k+1)}_2\\
0&\cdots &\cdots & 0& E^{(k+1)}_1
\end{pmatrix}.
\end{equation*}
This identity implies that $\psi$ induces a map of quotient rings
\begin{equation*}
\tilde\psi:\mathcal O(\Y_P)\To
\C\left[\,\si^{(1)}_{1},\dotsc,\si^{(k+1)}_{n-n_k}\, ,\,
q_{1},\dotsc,q_{k}\right]/(E^{(k+1)}_1,\dotsc, E^{(k+1)}_n).
\end{equation*}
This map together with \eqref{e:ON_Ppres} defines an inverse
$\mathcal N_P\to \Y_P$ to \eqref{e:Atilde}.
Thus $\tilde \psi$ is an isomorphism and everything follows from
Theorem \ref{t:Cio}.
\end{proof}

\begin{rem}
Actually, Peterson's results are more generally stated over the
integers. And here in particular the analogous theorem over $\Z$
holds, with the exact same proof. We have stayed over $\C$ in our
presentation since that is all we will require.
\end{rem}

\subsection{Toeplitz matrices}\label{s:Toeplitz}

The stabilizer of $f$ under conjugation by $U^-$ is precisely the
$(n-1)$-dimensional abelian subgroup of lower-triangular
unipotent Toeplitz matrices,
\begin{equation*}
 X:=(U^-)_f=\left\{x=\left .
\begin{pmatrix}
 1& & & & \\
 a_1& 1   &    &  &        \\
 a_2   & a_1 & \ddots&      &   \\
 \vdots   &     & \ddots & 1   &     \\
 a_{n-1}   &\cdots     &   a_2  & a_1   & 1
\end{pmatrix} \right | a_1,\dotsc, a_{n-1}\in \C\ \right\}.
\end{equation*}
Let us take the matrix entries $a_1,\dotsc, a_{n-1}$ as
coordinates on $ X$, thereby identifying  $\mathcal O(
X)=\C[a_1,\dotsc, a_{n-1}]$. For $1\le m\le n-1$ let
$\Delta_m\in\mathcal O( X)$ be defined by
\begin{equation}\label{e:dels}
\Delta_m=\det(a_{j-i+m})_{i,j=1}^{n-m}=\left|\begin{matrix} a_{m}
& a_{m-1}&\cdots &         \\
 a_{m+1}& \ddots& \ddots& \vdots \\
 \vdots   &\ddots  &\ddots &a_{m-1}\\
 a_{n-1}      & \cdots       &a_{m+1} &a_m\\
\end{matrix}\right |,
\end{equation}
where $a_0=1$ and $a_l=0$ if $l<0$.
Let
\begin{equation*}
 X_P=  X\cap B^+w_Pw_0 B^+.
\end{equation*}
We recall the following explicit description of the $ X_P\subset
 X$.
\begin{lem}[\cite{Rie:QCohGr} Lemma 2.5]
As a subset (not subvariety) of $ X$,  $ X_P$ is described by
\begin{equation*}
 X_P=\left \{\, u\in  X\ \left |\  \Delta_i(u)\ne 0\
\iff\ i\in\{n_1,\dotsc, n_k \}\right.\, \right\}.
\end{equation*} \qed
\end{lem}
Note that the map
\begin{equation*}
\Delta=(\Delta_1,\dotsc,\Delta_{n-1}):X\to\C^{n-1}
\end{equation*}
has the property $\Delta\inv(0)=\{\, 0\,\}$ (in fact over any
field). And $\Delta^*:\mathcal O(\C^{n-1})=\C[z_1,\dotsc,
z_{n-1}]\to \C[a_1,\dotsc , a_{n-1}]$ is homogeneous, if the
generators are taken with suitable degrees.
Therefore
$\Delta$ is a finite morphism, see e.g. \cite{GrHa:AlgGeom}.

\begin{thm}[D. Peterson]
\begin{enumerate}\item
Define
\begin{equation*}
\X:=\Y\cap\left( B^-w_0B^-/B^-\right) \quad\text{and}\quad \X_P:=
\X\cap \Y^P.
\end{equation*}
Then the isomorphism $U^-\to B^-w_0B^-/B^-$ defined by $u\mapsto
uw_0B^-$ identifies $X$ with $\X$ and also $X_P$ with $\X_P$
for each parabolic $P$.
\item
The map \eqref{e:PetIso1} induces an isomorphism of $\mathcal
O(\X_P)$ with $qH^*(G/P)[q_1\inv,\dotsc, q_k\inv]$ giving
\begin{equation}
\mathcal O(X_P)\overset\sim\longleftarrow\mathcal
O(\X_P)\overset{\sim}\To qH^*(G/P)[q_1\inv,\dotsc, q_k\inv].
\label{e:PetIso2}
\end{equation}
In particular, each $\X_P$ is
open dense in $\Y_P$.
\end{enumerate}
\end{thm}

For a proof of this when $P=B$ see Theorems~8 and 9 in
\cite{Kos:QCoh}. The general case is analogous, and for (2) see
also the proof of Lemma~\ref{l:GenKos} below.

\section{The Quantum parameters as functions on $X_P$}\label{s:Deltas}

After applying \eqref{e:PetIso2}, the quantum parameters $q^P_j$
may be expressed (up to taking some roots) in terms of the
functions $\Delta_i$ from \eqref{e:dels}. In the full flag
variety case, that is on $X_B$ where all $\Delta_i$ are
non-vanishing, Kostant \cite{Kos:QCoh} has given the following
formula,
\begin{equation}\label{e:Kos}
q_j^B=\frac{\Delta_{j-1}\Delta_{j+1}}{(\Delta_j)^2}.
\end{equation}
This generalizes as follows to the partial flag variety case.
\begin{lem}\label{l:GenKos}
Let the quantum parameters $q_j^P$ be regarded as functions on $
X_P$ via the isomorphism $ \mathcal O( X_P)\cong
qH^*(G/P)[q_1\inv,\dotsc,q_k\inv]$ from \eqref{e:PetIso2}. Then
\begin{equation*}
\left(q_j^P\right)^{(n_j-n_{j-1})(n_{j+1}-n_j)}=
\frac{(\Delta_{n_{j-1}})^{n_{j+1}-n_j}(\Delta_{n_{j+1}})^{n_{j}-n_{j-1}}}
{(\Delta_{n_j})^{n_{j+1}-n_{j-1}}}.
\end{equation*}
\end{lem}
The proof of this lemma which we give below is an adaptation of
Kostant's proof of the formula \eqref{e:Kos}.

\begin{proof}
Let $y\in X_P$. Then $yw_0B^-\in \mathcal X_P$ and we have
\begin{equation*}
yw_0B^-=uw_P B^-\quad \text{for some $u\in U^+$.}
\end{equation*}
Without loss of generality  $u$ may be chosen such that
$\Ad(u\inv)\cdot f$ is an ASK-matrix $A\in \mathcal N_P$. Since
$uw_PB^-=yw_0B^-$ we can find $\bar u\in U^+$ and $t\in T$ such
that
\begin{equation*}
y= uw_P\inv w_0\inv t \bar u\inv.
\end{equation*}
We have
\begin{equation}\label{e:demo}
\Ad (\bar u\inv)\cdot f=\Ad(t w_0w_P\inv  u\inv y)\cdot f= \Ad(t
w_0 w_P\inv )\cdot A
\end{equation}
The right hand side may be expanded to give
\begin{equation*}
\Ad(t\inv w_0w_P )\cdot A= \sum_{i=1}^{n-1} m_i f_i + {\text
{higher weight space terms,}}
\end{equation*}
where explicitly
\begin{equation*}
m_i=\begin{cases} (-1)^{n_{j+1}-n_j}\alpha_{i}(t)q^P_j(y) &\text{if } i=n-n_j\\
                  \alpha_{i}(t) &\text{if $n-i\notin \{n_1,\dotsc, n_k\}.$}
\end{cases}
\end{equation*}
This follows from the isomorphisms in Section~\ref{s:Peterson}. On
the other hand since $\bar u\in U^+$, the left hand side of
\eqref{e:demo} implies that all the $m_i$ must equal to $1$.
Therefore we have the identities
\begin{alignat*}{2}
\alpha_{n-n_j}(t)&= (-1)^{n_{j+1}-n_j}q^P_j(y)\inv & \qquad&
j=1,\dotsc, k, \\
\al_i(t)&=1 & &n-i\notin \{n_1,\dotsc , n_k\}.
\end{alignat*}
Thus $t$ is determined up to a scalar factor $\lambda$ by the
$q_i(y)$'s. Let $T_{k+1}$ be the $(n-n_k)\x (n-n_{k})$ identity
matrix and $T_j$ the $(n_{j}-n_{j-1})\x (n_{j}-n_{j-1})$ matrix
\begin{equation*}
T_j=(-1)^{n-n_j}q_{j}(y)q_{j+1}(y)\cdots q_k(y)
\begin{pmatrix}
1 &  &\\
               & \ddots  &\\
& & 1
\end{pmatrix},
\end{equation*}
for $1\le j\le k$. Then $t$ may explicitly be described by
\begin{equation*}
t=\lambda\begin{pmatrix}\ \fbox{\rule[-4mm]{0cm}{9mm} $T_{k+1}$}  & &\\
& \hskip -3.5mm \text{\fbox{\rule[-1.5mm]{0cm}{5mm} $T_{k}$}}    & &\\
& & \ddots & \\
& &        & \hskip -1.5mm\fbox{\rule[-1.5mm]{0cm}{5mm} $T_{1}$}\
\end{pmatrix}.
\end{equation*}
Now
\begin{multline*}
\Delta_{n_j}(y)=(y\cdot v_1\wedge\cdots \wedge v_{n-n_j+1}\ |\
v_{n_j+1}\wedge\cdots\wedge v_n)= \\ =(uw_P\inv w_0\inv t \bar
u\inv\cdot v_1\wedge\cdots \wedge v_{n-n_j}\ |\
v_{n_j+1}\wedge\cdots\wedge v_n)=\\ = \omega_{n-n_j}(t)(u w_P\inv
w_0\inv \cdot v_1\wedge\cdots \wedge v_{n-n_j}\ |\
v_{n_j+1}\wedge\cdots\wedge v_n)=\\=
\lambda^{n-n_j} q_{k}(y)^{n_k-n_j}q_{k-1}(y)^{n_{k-1}-n_j}\cdots
q_{j+1}(y)^{n_{j+1}-n_j}
\end{multline*}
and the identity follows.
\end{proof}

\section{Total Positivity}\label{s:totpos}

A matrix $A$ in $GL_n(\R)$ is called {\it totally positive} (or
{\it totally nonnegative}) if all the minors of $A$ are positive
(respectively nonnegative). In other words $A$ acts by positive or
nonnegative matrices in all the fundamental representations
$\bigwedge^k\R^n$ (with respect to their standard bases). These
matrices clearly form a semigroup. The concept of totally positive
matrices is in this sense more fundamental than the naive concept
simply of matrices with positive entries, which overemphasizes the
standard representation. Total positivity for $GL_n$ was mainly
studied in and around the 1950's by Schoenberg, Gantmacher-Krein,
Karlin and others, and has relationships with diverse applications
such as oscillating mechanical systems and planar Markov
processes.

More recently, G. Lusztig  \cite{Lus:TotPos94} extended the theory
of total positivity to all reductive algebraic groups. His
extension  rests around a beautiful connection with canonical
bases (for ADE type). This point of view on total positivity was
part of the motivation for the main result of this paper, stated
in the next section. It goes as follows.

Let us consider the lower uni-triangular matrices $U^-$ (staying
with $G=GL_n$, to avoid making further definitions). Let $U^-_{\ge
0}$ be the set of totally nonnegative matrices in $U^-$. And
define the `totally positive' part of $U^-$ by
\begin{equation*}
U^-_{> 0}:= U^-_{\ge 0}\cap B^+w_0 B^+.
\end{equation*}
Now the canonical basis of the quantized universal enveloping
algebra $\mathcal U_q^-$ defined by Lusztig and Kashiwara gives
rise, after dualizing and taking the classical limit, to a basis
$\mathcal B$ of the coordinate ring $\mathcal O(U^-)$. Lusztig
proved that the canonical basis has positive structure constants
for multiplication and comultiplication, using his geometric
construction of $\mathcal U_q^-$.
This is the main
ingredient for the following theorem.

\begin{thm}[\cite{Lus:TotPosCan}, see also Section 3.13 in
\cite{Lus:PartFlag}] Suppose $u\in U^-(\R)$. Then
\begin{equation*}
u\in U^-_{> 0}   \qquad \iff \qquad  b(u)> 0 \quad\text {for all
$b\in\mathcal B$.}
\end{equation*}
\end{thm}

The functions $b\in \mathcal B$ are matrix coefficients of $U^-$
in irreducible representations of $GL_n$ with respect to canonical
bases of these representations (obtained from the canonical basis
of $\mathcal U^-_q$). This theorem is a reformulation of a result
from \cite{Lus:TotPosCan}, which holds for any simply laced
reductive algebraic group.

Philosophically, $U^-$ is a variety with a special basis on its
coordinate ring (even a $\Z$-basis if we were to define $U^-$
over the integers) which moreover has nonnegative integer
structure constants. And by Lusztig's theorem, the question after
which $u\in U^-$ have the property that all $b\in\mathcal B$ are
positive on $u$ has a very nice answer, namely the totally
positive part of $U^-$.

We ask the same question for the components $\mathcal Y_P$ of the
Peterson variety (or equivalently for the $X_P\subset U^-$),
whose coordinate rings are naturally endowed with Schubert bases
also with positive structure constants (as enumerative
Gromov-Witten invariants). And remarkably we discover total
positivity again in the answer.

The corollary, the parameterization result for totally nonnegative
finite Toeplitz matrices stated in the introduction, also
illustrates a common feature in total positivity. For example
there are natural parameterizations of $U^-_{>0}$ (introduced in
\cite{Lus:TotPos94}) which are related to combinatorics of the
canonical basis and have been studied extensively. See e.g.
\cite{FoZe:Intelligencer} for a survey.

\section{Statement of the main theorem}\label{s:main}

The varieties we are studying lie either inside $GL_n$ or
$G/B^-$. By their real points we mean coming from their split
real form, $GL_n(\R)$ and the real flag variety. We consider the
real points always to be endowed with the usual Hausdorff topology
coming from $\R$. The positive parts will be semi-algebraic
subsets of the real points.

Following \cite{Lus:TotPos94}, the totally positive part
$(G/B^-)_{>0}$ of $G/B^-$ is defined as the image of $U^+_{>0}$
under the quotient map $G\to G/B^-$. By a result in
\cite{Lus:TotPos94} this agrees with the image of
$U^-_{>0}w_0B^-$. So
\begin{equation*}
(G/B^-)_{>0}=U^+_{>0} B^-/B^-=U^-_{>0} w_0 B^-/B^-.
\end{equation*}
The totally nonnegative part $(G/B^-)_{\ge 0}$ is the closure of
$(G/B^-)_{>0}$ inside the real flag variety.

Using Peterson's isomorphisms \eqref{e:PetIso1} and
\eqref{e:PetIso2} we may evaluate elements of $qH^*(G/P)$ as
functions on the points of $\Y_P$ and $\mathcal X_P$, or $X_P$.
\begin{defn}
Let the {\it totally positive} part of $\Y_P$ be defined as
\begin{equation*}
\Y_{P,>0}:= \Y_P(\R)\cap (G/B^-)_{>0}.
\end{equation*}
This automatically lies in $\X_P$, so we also set
$\X_{P,>0}:=\Y_{P,>0}$. Finally, compatible with this definition,
set $X_{P,>0}:=X_P(\R)\cap U^-_{>0}$.

We define the {\it Schubert-positive} parts of $\Y_P, \X_P$ and
$X_P$, also compatibly with the various morphisms between them, by
\begin{align*}
\Y_{P,>0}^{Schub}&:=\{x\in \Y_{P}(\R)\ |\ \si_w^P(x)>0\ \text{ all
}w\in W^P \},\\
\X_{P,>0}^{Schub}&:=\{x\in \X_{P}(\R)\ |\ \si_w^P(x)>0\ \text{ all
}w\in W^P \},\\
X_{P,>0}^{Schub}&:=\{x\in X_{P}(\R)\ |\ \si_w^P(x)>0\ \text{ all
}w\in W^P \}.
\end{align*}
These are all semi-algebraic subsets of the real points of $\Y_P$,
$\X_P$ and $X_P$, respectively.
\end{defn}

\begin{thm}\label{t:main}
\begin{enumerate}
\item
The ramified cover $\pi=\pi^P=(q^P_1,\dotsc, q^P_k):\Y_P(\C) \to
\C^k$ restricts to a bijection
\begin{equation*}
\pi^P_{>0}:\Y^{Schub}_{P,>0}\to\R^k_{>0}.
\end{equation*}
\item
$\Y^{Schub}_{P,>0}$ lies in the smooth locus of $\Y_{P}$, and the
inverse of the map $\pi^P_{>0}:~\Y^{Schub}_{P,>0}\to\R^{k}_{>0}$
is analytic.
\item
The two notions of positivity agree. That is,
\begin{equation*}
\Y^{Schub}_{P,>0}=\Y_{P,>0},
\end{equation*}
and also $\Y^{Schub}_{P,>0}=\X^{Schub}_{P,>0}=\X_{P,>0}$ and
$X^{Schub}_{P,>0}=X_{P,>0}$.
\end{enumerate}
\end{thm}

\begin{rem}
Bianalyticity of $\pi^B_{>0}$ is equivalent to the non-vanishing
on $X_{B,>0}$ of the {\it quantum Vandermonde} function defined by
Kostant \cite{Kos:QCoh}, Section~9. This function on $X_B$ is
expressed as determinant of a matrix whose entries are
alternating sums of minors. By Theorem~\ref{t:main}.(2) it must
take either solely positive or negative values on $X_{B,>0}$. We
expect that the values will always be positive, but have checked
this so far only in very low rank cases.

We conjecture that all the analogous results to those stated in
Theorem~\ref{t:main} should hold true in general type. The
stabilizer of the principal nilpotent $f$ in that case should
have a totally nonnegative part with a cell decomposition coming
from Bruhat decomposition. And there should be an analogous
relationship with Schubert bases for quantum cohomology rings
$qH^*(G^\vee/P^\vee)$ of partial flag varieties of the Langlands
dual group, via the Peterson variety for general type.
\end{rem}

We now deduce the corollary stated as Proposition~\ref{t:2} in
the introduction.
\begin{cor}
Let $X_{\ge 0}$ denote the semi-algebraic subset of $X(\R)$ of
totally nonnegative unipotent lower-triangular Toeplitz matrices.
Then the restriction
\begin{equation*}
\Delta_{\ge 0}: X_{\ge 0}\To \R_{\ge 0}^{n-1}
\end{equation*}
of $\Delta:=(\Delta_1,\dotsc, \Delta_{n-1})$ is a homeomorphism.
\end{cor}

\begin{proof}[Proof of the Corollary]
By Theorem \ref{t:main} we have homeomorphisms
\begin{equation*}
(q^P_1,\dotsc, q^P_k):X_{P,>0}\to \R^{k}_{>0},
\end{equation*}
one for each parabolic $P$. By Lemma \ref{l:GenKos} the $q_j^P$'s
are related to the $\Delta_{n_j}$'s by a transformation which is
continuously invertible over $\R_{>0}^{k}$. By this observation,
and since $X_{\ge 0}=\bigsqcup X_{P,>0}$, we have that
\begin{equation*}
\Delta_{\ge 0}: X_{\ge 0}\To \R_{\ge 0}^{n-1}
\end{equation*}
is bijective. So $\Delta_{\ge 0}$ is continuous, bijective, and a
homeomorphism onto its image when restricted to any $X_{P,>0}$.
Since $\Delta$ is finite it follows that $\Delta_{\ge 0}\inv$ is
also continuous.
\end{proof}

\section{Grassmannians etc.}\label{s:Grassmannian}
The quantum cohomology rings of Grassmannians have been studied
much more extensively than those of partial flag varieties, see
for example \cite{Be:QSchubCalc,Gepner:FusRing,SiTi:QCoh,
Witten:VerAlgGras}. The Grassmannian case can be considered as a
kind of toy model for this paper. The main results stated in the
previous section are generalizing properties from the Grassmannian
case which were studied by elementary means, basically playing
with Schur polynomials, in \cite{Rie:QCohGr}. We will briefly
recall what happens in that case.

\begin{enumerate}
\item
Let $\mathcal V_{d,n}$ be the transpose  of $\bar X_{P_d}$
(upper-triangular rather than lower-triangular Toeplitz matrices),
notation as in \cite{Rie:QCohGr}. And let $\mathcal
O_{red}(\mathcal V_{d,n})$ be the reduced coordinate ring of
$\mathcal V_{d,n}$. Then an incarnation of Peterson's theorem says
$\mathcal O_{red}(\mathcal V_{d,n}) \cong qH^*(G/P_d)$.
\item
The points of $\mathcal V_{d,n}$ are those
\begin{equation*}
u=\begin{pmatrix}
 1 & a_1&\cdots & a_d   &        &0 \\
   &1   & a_1   &       &\ddots  & \\
   &    &  \ddots &\ddots &      & a_d \\
   &     &        & \ddots & a_1 &\vdots  \\
   &     &        &      &  1    &a_1 \\
   &     &       &       &       &1
\end{pmatrix}
\end{equation*}
for which
\begin{equation*}
p(x)=x^d+ a_1 x^{d-1}+\dotsc + a_d=\prod_{j=1}^d\left (x + z
e^{\left ( m_j\frac{2\pi i }{n}\right )} \right )
\end{equation*}
for some $z\in\C$ and integers $0\le m_1<\dotsc< m_d< n$. In
other words either $u$ is the identity matrix or otherwise the
roots of the generating polynomial $p(x)$ are distinct complex
numbers with $x_1,\dotsc, x_d$ with $x_1^n=\dotsc=x_d^n$. Write
$u=u(x_1,\dotsc, x_d)$.
\item
Let $u=u(x_1,\dotsc, x_d)$ as above, and $w_\lambda$ the
Grassmannian permutation in $W^{P_d}$ corresponding to a Young
diagram $\lambda$. The image of the Schubert class
$\si^{P_d}_{w_\la}\in \mathcal O_{red}(\mathcal V_{d,n})$ is given
by
\begin{equation*}
\si^{P_d}_{w_\la}(u)=s_\lambda(x_1,\dotsc, x_d),
\end{equation*}
where $s_\lambda$ is the Schur polynomial associated to $\lambda$.
\item
Let $\zeta=e^{\frac {2\pi i}n}$, and set $u_{\ge 0
}(t)=u(t\zeta^{-\frac{d-1}2},t\zeta^{-\frac{d-1}2+1},\dotsc,t\zeta^{\frac{d-1}2})$.
Then
\begin{equation*}
u_{\ge 0}:\R_{\ge 0}\overset\sim\To  \left(\mathcal
V_{d,n}\right)_{\ge 0}\quad : \quad t\mapsto u_{\ge 0}(t)
\end{equation*}
is a homeomorphism, where $\left(\mathcal V_{d,n}\right)_{\ge 0}$
denotes the totally nonnegative matrices in $\mathcal V_{d,n}$.
\item
The values of the Schubert classes on the $u(t)$ are given by a
closed (hook-length) formula, which explicitly shows them to be
positive for $t>0$.
\item The quantum parameter $q$ is given by
$q(u(x_1,\dotsc,x_n))=(-1)^{d+1}x_1^n$. In particular, $q(u_{\ge 0
}(t))=t^n$.
\end{enumerate}

From the proof of Theorem~\ref{t:Pet}, in particular from
inspection of the matrix $u$ introduced in \eqref{e:u}, we see
directly that $\mathcal V_{d,n}\cong\Y_{P_{d}}$ via $u\mapsto
uw_{P_d} B^-$. Notice also that (4) and (6) give the
parameterization by quantum parameters of Theorem~\ref{t:main} in
this special case.

Peterson has announced in \cite{Pet:QCoh} that all the quantum
cohomology rings $qH^*(G/P)$ are reduced. To prove this amounts to
showing that the element $\sum_{w\in W^P}\sigma_w \sigma_{PD(w)}$
is a nonzerodivisor in $qH^*(G/P)$ (see also
\cite{Abrams:QEuler}). This is because, for example, if $\sigma\in
qH^*(G/P)$ is nilpotent then all $\mu\si$ for $\mu\in qH^*(G/P)$
are, and the corresponding multiplication operators $M_{\mu\si}$
on $qH^*(G/P)$ have vanishing trace. But computing these traces by
Poincar\'e duality gives $\tr(M_{\mu\si })=\left<\mu\si\sum_{w\in
W^P} \si_w\, \si_{PD(w)}\right>_{\mathbf q}=0$ and therefore
$\si\cdot\left( \sum_{w\in W^P }\si_w\si_{PD(w)}\right)=0$. It is
in fact sufficient to show that $\sum_{w\in W^P} \si_w
\si_{PD(w)}$ is generically nonvanishing on $\mathcal Y_P$. See
e.g. Lemma~\ref{l:Jacobian} or the direct linear algebra proof
from \cite{Rie:QCohGr}~Section~5.2.

Apart from the Grassmannian case where everything is very
explicit, and the full flag variety case treated in
\cite{Kos:QCoh}, where the Peterson variety $\mathcal Y_B$ is
irreducible, I do not know a proof that $\sum_{w\in W^P} \si_w
\si_{PD(w)}$ is generically nonvanishing. But with the help of the
explicit results above we can prove at least the following lemma
which will come in handy later.

\begin{lem}\label{l:topclass}
The element of $\mathcal O(\Y_P)$ defined by $\si^P_{w_0^P}$
takes nonzero values on an open dense subset of $\Y_P(\C)$.
\end{lem}
\begin{proof}
Since $\X_P$ is open dense in $\Y_{P}$ it suffices to show that
$\si^P_{w_0^P}$ is nonzero on an open dense subset of $\X_P(\C)$.
By Theorem~\ref{t:Pet}(3) we may furthermore replace $\X_P$ by
$X_P$. So let us identify the Schubert classes $\sigma^P_w$ with
rational functions on $\bar X_P$ and prove that the top one is
generically non-vanishing.

Let $P_{m}$ denote the maximal parabolic with $I^{P_m}=\{m\}$ and
let $C$ be an irreducible component of the closure $\bar
X_P=\bigsqcup_{P'\supseteq P}X_{P'}$. If $I^P=\{i_1,\dotsc,
i_k\}$ then we have
\begin{equation*}
(\Delta_{n_1},\dotsc,\Delta_{n_k}): \bar X_P(\C)\To \C^k
\end{equation*}
is finite, as pullback of the finite map
$(\Delta_{1},\dotsc,\Delta_{n-1}): X(\C)\To \C^{n-1}$. Therefore
the restriction of $(\Delta_{n_1},\dotsc, \Delta_{n_k})$ to $C$ is
surjective and $C$ intersects all of the subvarieties $\bar
X_{P_{n_i}}$ of $\bar X_P$.

Now in $qH^*(G/P)$ we have
\begin{equation*}
\si^P_{w^P}=\si^P_{s_1\cdots s_{n_1}}\cdot\si^P_{s_1\cdots
s_{n_2}}\cdot\dotsc\cdot \si^P_{s_1\cdots s_{n_k}}.
\end{equation*}
Let $x\in X_P$. Then tracing through Peterson's isomorphisms gives
\begin{equation*}
\si^P_{s_1\cdots s_{n_j}}(x)=G^{n_j}_{n_j}(x w_0 B^-),
\end{equation*}
where $G^{n_j}_{n_j}$ is as in \eqref{e:Gs}. This function extends
to $X_{P_{n_j}}\subset \bar X_P$ and is seen to be non-vanishing
there using the explicit description of $X_{P_{n_j}}$ (see (2)
above). Since any irreducible component of $\bar X_P$ meets
$X_{P_{n_j}}$, we have that $\si^P_{s_1\cdots s_{n_j}}$ is
generically nonzero on $X_P$. The same holds therefore for
$\si^P_{w^P}$ as the product of the $\si^P_{s_1\cdots s_{n_j}}$.
\end{proof}

\section{Proof of Theorem \ref{t:main}.(1)}\label{s:proof1}
We must first check that $q^P_{>0}$ actually takes values in
$\R_{>0}^{k}$. This follows from the following observation.
\begin{lem}\label{l:qpos}
Let $P_{n_j}$ be the maximal parabolic defined by
$I^{P_{n_j}}=\{n_j\}$, and set $v=w_0^{P_{n_j}}\in W^{P_{n_j}}$
to be the longest element. Then $v\in W^P$ and we have the
following relation in $qH^*(G/P)$,
\begin{equation*}
\si^P_{s_{n_j}}\cdot\si^P_{v}=q^P_j\ \si^P_{v\tau_{j,j}},
\end{equation*}
where $\tau_{j,j}=s_{n_j}\cdots s_{n_{j+1}-1} s_{n_j-1}\cdots
s_{n_{j-1}+1} $.
\end{lem}
\begin{proof}
Let $\alpha$ be a positive root such that
$<\alpha,\omega_{n_j}^\vee>\ne 0$. So $\alpha=\alpha_h+\dotsc +
\alpha_l$ for some  $h\le n_j\le l$. By the Chevalley formula,
$\sigma_{vs_\alpha}$ appears in the expansion of the product only
if $\ell(vs_\alpha)=\ell(v)+1$. If $h<n_j<l$ then
$\ell(vs_\alpha)=\ell(v)+\ell(s_\alpha)\ge \ell(v)+3$. So assume
$h=i$ or $l=i$. In either of those two cases
$\ell(vs_\alpha)=\ell(v)-1$. So the classical contribution to
$\si^P_{s_{n_j}}\cdot\si^P_{v}$ is indeed zero.

Suppose now $\ell(v\tau_{h,l})=\ell(v)-\ell(\tau_{h,l})$, where
$\tau_{h,l}$ is as in Section \ref{s:QChev}. This is equivalent to
asking $\tau_{h,l}\inv\in W^{P_{n_j}}$. Since
\begin{equation*}
\tau_{h,l}\inv=s_{n_{h-1}+1}\cdots s_{n_{l}-1}s_{n_{l+1}-1}\cdots
s_{n_h}
\end{equation*}
sends both $\alpha_{n_h}$ and $\alpha_{n_l}=s_{n_h}\cdots
s_{n_{l+1}-1}(\alpha_{n_l-1})$ to negative roots we must have
$h=l=j$. So by quantum Chevalley's rule the only possible quantum
contribution to the product $\si^P_{s_{n_j}}\cdot\si^P_{v}$ is
$q^P_j\ \si^P_{v\tau_{j,j}}$. It follows by a direct check that
this term does indeed appear (as of course it must, since the
product cannot be zero by the same arguments as in
Lemma~\ref{l:topclass}.)
\end{proof}
Now we would like to show that $\pi^P_{>0}$ is actually
surjective. For this fix a point $Q\in (\R_{>0})^{k}$ and
consider its fiber under $\pi=\pi^P$. We may regard
\begin{equation*}
R_Q:=qH^*(G/P)/(q^P_1-Q_1,\dotsc, q_k^P-Q_k)
\end{equation*}
as the (possibly non-reduced) coordinate ring of $\pi\inv(Q)$.
Note that $R_Q$ is a finite-dimensional algebra with basis given
by the (image of the) Schubert basis. We will use the same
notation $\si^P_w$ for the restriction of a Schubert basis element
to $R_Q$.
\begin{lem}\label{l:EV}
Suppose $\mu\in R_Q$ is a nonzero simultaneous eigenvector for all
linear operators $R_Q\to R_Q$ which are defined by multiplication by
elements in $R_Q$. Then there exists a point $p\in\pi\inv(Q)$ such that (up to
a scalar factor)
\begin{equation*}
\mu=\sum_{w\in W^P}\si^P_w(p)\,\si^P_{PD(w)}.
\end{equation*}
\end{lem}
\begin{proof}
Consider the algebra homomorphism
\begin{equation*}
R_Q\To \C
\end{equation*}
which takes $\si\in R_Q$ to its eigenvalue on the eigenvector
$\mu$. This defines the $\C$-valued point $p$ in $\pi\inv(Q)$.
Now let us write $\mu$ in the Schubert basis,
\begin{equation*}
\mu=\sum_{w\in W^P}m_w\si^{PD(w)}, \quad\qquad m_w\in\C.
\end{equation*}
For $\si\in R_Q$, let $\left<\sigma\right>_Q\in \C$ denote the coefficient of
$\sigma^P_{w_0^P}$ in the Schubert basis expansion of $\si$. Then
by quantum Poincar\'e duality we have
\begin{equation*}
m_w=\left<\,\si^w\cdot \mu\,\right >_Q=\left <\,\si^w(p)\
\mu\,\right>_Q=\si^w(p)\,\left<\,\mu\,\right>_Q=
\si^w(p)\, m_1.
\end{equation*}
Here $m_1$ must be a nonzero scalar factor (since $\mu\ne 0$), and the lemma is proved.
\end{proof}

We continue the Proof of Theorem \ref{t:main}.(1) our immediate
aim being to find a Schubert positive point $p_0$ in the fiber $\pi\inv(Q)$.
Set
\begin{equation*}
\sigma:=\sum_{w\in W^P}\si_w^P\in R_Q.
\end{equation*}
Suppose the multiplication operator on $R_Q$ defined by
multiplication by $\si$ is given by the matrix
$M_\si=(m_{v,w})_{v,w\in W^P}$ with respect to the Schubert
basis. That is,
\begin{equation*}
\si\cdot\si_v^P=\sum_{w\in W^P} m_{v,w}\si^P_w.
\end{equation*}
Then since $Q\in \R_{>0}^k$ and by positivity of the structure
constants it follows that $M_\si$ is a nonnegative matrix.
Furthermore let us assume the following lemma (to be proved
later).

\begin{lem}\label{l:indec}
$M_\si$ is an indecomposable matrix.
\end{lem}

Given the indecomposable nonnegative matrix $M_\si $,
then by Perron-Frobenius theory (see e.g. \cite{Minc:NonnegMat}
Section 1.4) we know the following. \vskip .3cm
\parbox[c]{12cm}{
The matrix $M_\si$ has a positive eigenvector $\mu$ which is
unique up to scalar (positive meaning it has positive coefficients
with respect to the standard basis). Its eigenvalue, called the
Perron-Frobenius eigenvalue, is positive, has maximal absolute
value among all eigenvalues of $M_{\sigma}$, and has algebraic
multiplicity $1$. The eigenvector $\mu$ is unique even in the
stronger sense that any nonnegative eigenvector of $M_\sigma$ is a
multiple of $\mu$.} \vskip .3cm

Suppose $\mu$ is this eigenvector chosen normalized such that
$\left<\mu\right>_Q=1$. Then since the eigenspace containing $\mu$
is $1$--dimensional, it follows that $\mu$ is joint eigenvector
for all multiplication operators of $R_Q$. Therefore by
Lemma~\ref{l:EV} there exists a $p_0\in \pi\inv(Q)$ such that
\begin{equation*}
\mu=\sum_{w\in W^P}\si^P_w(p_0)\, \si^P_{PD(w)}.
\end{equation*}
Positivity of $\mu$ implies that $\si^P_w(p_0)\in\R_{>0}$ for all
$w\in W^P$. Hence $p_0\in \Y_{P,>0}^{Schub}$. Also the point $p_0$
in the fiber with this property is unique.

Therefore we have shown modulo the Lemma~\ref{l:indec} that
\begin{equation}\label{e:homeo}
\Y^{Schub}_{P,>0}\To \R_{>0}^{k}
\end{equation}
is a bijection.
Finally we complete the proof of Theorem \ref{t:main}.(1) by
proving the lemma.
\begin{proof}[Proof of Lemma~\ref{l:indec}] Recall that $\si=\sum_{w\in
W^P}\si_w$. Suppose indirectly that the matrix $M_\si$ is
reducible. Then there exists a nonempty, proper subset $V\subset
W^P$ such that the span of $\{\si_v\ | \ v\in V \}$ in $R_Q$ is
invariant under $M_\si$. We will derive a contradiction to this
statement.

First take any element $v\in V$. Then the top class $\si_{w_0^P}$
occurs in $\si\cdot \si_v$ with coefficient $1$ by quantum
Poincar\'e duality. Therefore we have $w_0^P\in V$.

Next we deduce that $1\in V$. Suppose not. Then the coefficient of
$\si_{1}$ in $\si_{w}\cdot\si_{w_0^P}$ must be zero for all $w\in
W^P$, or equivalently
\begin{equation*}
\left<\,\si_w\cdot\si_{w_0^P}\cdot \si_{w_0^P}\,\right>_Q=0
\end{equation*}
for all $w\in W^P$. But this also implies
$\left<\,\si_w\cdot\si_{w_0^P}\cdot \si_{w_0^P}\,\right>_\mathbf
q=0$, since the latter is a nonnegative polynomial in the
$q_i^P$'s which evaluated at $Q\in \R_{>0}^k$ equals $0$.
Therefore $\si_{w_0^P}\cdot\si_{w_0^P}=0$ in $qH^*(G/P)$ by
quantum Poincar\'e duality. This leads to a contradiction with
Lemma~\ref{l:topclass}, that the element $\si_{w_0^P}$ is
generically nonzero as function on $\Y_P$.
\begin{footnote}{Fulton and Woodward \cite{FuWo:SchubProds} have in fact
recently proved that no two Schubert classes in $qH^*(G/P)$ ever
multiply to zero. This result can also be recovered as a corollary
of Theorem~\ref{t:main}, since the product of two Schubert classes
must take positive values on $\Y_{P,>0}$ and hence cannot be zero.
}
\end{footnote}

So $V$ must contain $1$. Since $V$ is a proper subset of $W^P$ we
can find some $w\notin V$. In particular, $w\ne 1$. It is a
straightforward exercise that given $1\ne w\in W^P$ there exists
$\alpha\in \Delta_+^P$ and $v\in W^P$ such that
\begin{equation*}
w=v s_\alpha, \quad\text{and}\quad \ell(w)=\ell(v)+1.
\end{equation*}
Now $\alpha\in \Delta_+^P$ means there exists $n_j\in I^P$ such
that $<\alpha,\omega_{n_j}^\vee>\ne 0$. And hence by the
(classical) Chevalley Formula we have that
$\si_{s_{n_j}}\cdot\si_v$ has $\si_w$ as a summand. But if
$w\notin V$ this implies that also $v\notin V$, since $\si\cdot
\si_v$ would have summand $\si_{s_{n_j}}\cdot\si_v$ which has
summand $\si_w$. Note that there are no cancellations with other
terms by positivity of the structure constants.

By this process we can find ever smaller elements of $W^P$ which
do not lie in $V$ until we end up with the identity element, so a
contradiction.
\end{proof}

\section{Proof of Theorem \ref{t:main}.(2)}\label{s:proof2}

We need to show that $\Y_{P,>0}$ lies in the smooth locus of $\Y_P$.
Consider the map
\begin{equation*}
E=\left [(E^{(k+1)}_1,\dotsc, E^{(k+1)}_n)\right ]:\ \C^n\To
\C[q_1,\dotsc, q_k]^n,
\end{equation*}
and its evaluation at $Q=(Q_1,\dotsc,Q_k)\in \C^k$,
\begin{equation*}
E_Q=\operatorname{ev}_{Q}\circ E
:\ \C^n\To \C^n.
\end{equation*}
Here the coordinates $\epsilon_1,\dotsc,\epsilon_n$ of the source
$\C^n$ are the $\sigma^{(m)}_i=:\epsilon_{n_{m-1}+i}$. Let
\begin{equation*}
J_E:=\det\left(\frac{\partial E^{(k+1)}_i}{\partial
\epsilon_j}\right)_{i,j}\in \C[\si^{(1)}_1,\dotsc,
\si^{(k+1)}_{n-n_k}\, ,\, q_1,\dotsc, q_k],
\end{equation*}
which at $\mathbf q=Q$ evaluates to
$J_{E_Q}=\det\left(\frac{\partial (E_{Q})_i}{\partial
\epsilon_j}\right)_{i,j}\in\C[\sigma^{(1)}_1,\dotsc,
\sigma^{(k+1)}_{n-n_k}]$, the Jacobian of $E_Q$. Let us also
denote by $J_{E}$ and $J_{E_Q}$ the classes these functions
define via \eqref{e:presentation} in $qH^*(G/P)$ and in
$R_Q=qH^*(G/P)/(q_1-Q_1,\dotsc, q_k-Q_k)$, respectively.

Note that the zero-fiber of $E_Q$ equals $(\pi^P)\inv(Q)$, and a
point $p\in (\pi^P)\inv(Q)$ is a smooth point of $\Y_P$ if the
$J_{E_Q}(p)\ne 0$. The smoothness assertion of
Theorem~\ref{t:main}.(2) follows from the following lemma.

\begin{lem}\label{l:Jacobian} The element
$J_E\in qH^*(G/P)$ is expressed in terms of the Schubert basis by
\begin{equation}\label{e:Jacobian}
J_{E}\ =\ \sum_{w\in W^P} \si_w \si_{PD(w)}.
\end{equation}
\end{lem}

\begin{proof}[Proof]
The main ingredient for this lemma is a result from
\cite{CatDickSturm:residues} or \cite{SiTi:QCoh}. But we begin by
checking the normalization. Following \cite{Kim:QCohPFl} we have
$\left <J_E\right>_\mathbf q=\left <J_{E_Q}\right>_Q=|W^P|$. In
fact, in terms of the Chern roots $J_{E_0}$ is expressed
explicitly by
\begin{equation*}
J_{E_0}=\prod_{\begin{smallmatrix} (i,j), \ \text{s.t.}\  i\le n_m <j \\
\text{some $1\le m\le k$ }\end{smallmatrix}}(x_i- x_j)\in
\C[x_1,\dotsc, x_{n}]^{W_P}\cong\C[\sigma^{(1)}_1,\dotsc,
\sigma^{(k+1)}_{n-n_k}],
\end{equation*}
and hence represents the Euler class in $H^*(G/P)$. Therefore,
$\left <J_{E_{0}}\right >_0=\int_{G/P}\chi_{G/P}=| W^P |$. But by
its degree $\left <J_E\right>_\mathbf q=\left <J_{E_Q}\right>_Q$
is a constant, independent of $Q$.


Now given the normalization as above, \cite{SiTi:QCoh}
Proposition~4.1 says that
$\operatorname{Res_{E_Q}}(\tilde\eta)=\left<\eta \right>_Q$,
where $\tilde\eta\in\mathcal O(\C^n)$ and $\eta\in
R_Q\cong\mathcal O(\C^n)/((E_Q)_1,\dotsc,(E_Q)_n)$ is the class
represented by $\tilde\eta$. Putting this identity together with
\cite{SiTi:QCoh} Lemma~4.3 we obtain the identity
\begin{equation*}
\tr(M_\kappa) = \left<\kappa
J_{E_{Q}}\right>_Q,\qquad\text{$\kappa\in R_Q$},
\end{equation*}
where $M_\kappa$ is the multiplication operator by $\kappa$ on
$R_Q$.

On the other hand this trace may be computed
from Poincar\'e duality by
\begin{equation*}
\tr(M_\kappa)=\left<\kappa\sum_{w\in W^P} \si_w
\si_{PD(w)}\right>_Q.
\end{equation*}
Comparing the two expressions for all $Q$ and all $\kappa$ it
follows that
\begin{equation}
J_{E}\ =\ \sum_{w\in W^P} \si_w \si_{PD(w)}
\end{equation}
as required.
\end{proof}

\vskip .2cm It remains to prove that the inverse to $\pi^P_{>0} $
is analytic. This follows from the following lemma.

\begin{lem}\label{l:analityc} Choose local coordinates $y_1,\dotsc, y_k$ in a
neighborhood of $p_0\in X_{P,>0}$. The Jacobian $\mathcal J=\det
\left (\frac{\partial q^P_i}{\partial y_j}\right)$ is nonzero at
the point $p_0$.
\end{lem}

\begin{proof}
Let $Q=\pi^P(p_0)$. Let $R=qH^*(G/P)$ and $I\subset R$ the ideal
$(q_1-Q_1,\dotsc, q_k-Q_k)$.
The Artinian ring $R_Q=R/I$ is isomorphic to the sum of local
rings $R_Q\cong\bigoplus_{x\in (\pi^P)\inv(Q)} R_x/IR_x$. And for
$x=p_0$ the local ring $R_{p_0}/IR_{p_0}$ corresponds in $R_Q$ to
the Perron--Frobenius eigenspace of the multiplication operator
$M_\sigma$ from the above proof. Since this is a one-dimensional
eigenspace (with algebraic multiplicity one) we have that
$\dim(R_{p_0}/IR_{p_0})=1$. Therefore any non-zero element $r\in
R_{p_0}/IR_{p_0}$ has the property $r(p_0)\ne 0$. But the Jacobian
$\mathcal J$  gives a non-trivial element in $R_{p_0}/IR_{p_0}$ since
its residue at $p_0$ with respect to $I$ is nonzero 
(see e.g. Chapter 5 in \cite{GrHa:AlgGeom}). 
\end{proof}

\section{The Schubert classes as rational
functions on $\mathcal Y$}\label{s:SchubertClasses}

To compare Schubert-positivity with total positivity we need to
make a closer study of the functions defined by the Schubert
classes. The following proposition is one of the most striking
features of the Peterson variety picture of quantum cohomology.
As far as I understand, it can be extracted from Peterson's
statements in \cite{Pet:QCoh} or \cite{Pet:Montreal} on the
connection between each of the $qH^*(G/P)$'s and the homology of
the loop group $\Omega K$ of the compact real form of $G$. We will
give a direct proof here for type $A$.

\begin{prop}[D. Peterson]\label{p:SchubRestr}
Let $w\in W$ and $\si^B_w$ the corresponding Schubert class
considered as a function on $\Y_B$. Let $\widetilde{\si}_w$ be
the rational function on the Peterson variety $\Y=\overline\Y_B$
defined by $\widetilde{\si}_w|_{\Y_B}=\si^B_w$. Then
$\widetilde{\si}_w$ is regular on $\Y_P\subset \Y$ if
$w\in W^P$. And in that case we have
\begin{equation*}
\widetilde{\si}_w|_{\Y_P}=\si^P_w\ \in \ \mathcal O(\Y_P).
\end{equation*}
\end{prop}
Our proof of this proposition uses the following lemma.

\begin{lem}\label{l:vanishing}  Suppose that $j\in I^P$ and $j+1\notin
 I^P$. Then the rational function $q^B_j\left(G^{j}_{i-1} G^{j-1}_{l-1}-
G^{j-1}_{i-2} G^{j}_{l}\right)$ vanishes on $\Y_P$.
\end{lem}
\begin{proof}
Let $gB^-\in\Y_P$. Then $(g\cdot v_{-\omega_m}\ |\
v_{-\omega_m})\ne 0$ precisely if $m\in I^P$, and in this case
\begin{equation*}
G^m_{i}(gB^-)=\frac{(g\cdot v_{-\omega_m}\ |\ s_{m-i+1}\cdots
s_{m}\cdot v_{-\omega_m}) } {(g\cdot v_{-\omega_m}\ |\
v_{-\omega_m})}
\end{equation*}
is well defined.  Also \eqref{e:Kos} implies that $q^B_m$ is well
defined on $\Y_P$ whenever $m\in I^P$, and is given by
\begin{equation*}
q^B_m(gB^-)=\frac{(g\cdot v_{-\omega_{m-1}}\ |\ v_{-\omega_{m-1}})
(g\cdot v_{-\omega_{m+1}}\ |\ v_{-\omega_{m+1}}) } {(g\cdot
v_{-\omega_m}\ |\ v_{-\omega_m})^2}.
\end{equation*}
Therefore we have
\begin{multline*}
q^B_j\left( G^{j}_{i-1} G^{j-1}_{l-1}- G^{j-1}_{i-2}
G^{j}_{l}\right)= \frac {(g\cdot v_{-\omega_{j+1}}\ |\
v_{-\omega_{j+1}} )}{ (g\cdot v_{-\omega_j}\ |\ v_{-\omega_j})^3}\cdot\\
\cdot\big ( (g\cdot v_{-\omega_j}\ |\ s_{j-i+2}\cdots s_{j}\cdot
v_{-\omega_j}) (g\cdot v_{-\omega_{j-1}}\ ,\ s_{j-l+1}\cdots
s_{j-1}\cdot v_{-\omega_{j-1}})- \\ \qquad\qquad - (g\cdot
v_{-\omega_j}\ |\ s_{j-i+2}\cdots s_{j-1}\cdot v_{-\omega_{j-1}})
(g\cdot v_{-\omega_{j-1}}\ |\ s_{j-l+1}\cdots s_{j}\cdot
v_{-\omega_{j}}) \big )
\end{multline*}
Now $(j+1)\notin I^P$ and $j\in I^P$ implies that $(g\cdot
v_{-\omega_{j+1}}\ |\ v_{-\omega_{j+1}})=0$  while $(g\cdot
v_{-\omega_j}\ |\ v_{-\omega_j})\ne 0$ on $\Y_P$. Hence the above
expression vanishes on $\Y_P$.
\end{proof}

\begin{proof}[Proof of Proposition \ref{p:SchubRestr}]
If $w=s_{h-i+1} \cdots s_{h-1}s_{h}$, then we have
\begin{equation*}
\widetilde \si_{s_{h-i+1}\cdots s_h}=G_{s_{h-i+1}\cdots s_h}=G^h_i
\end{equation*}
and the Proposition holds in this case by Theorem~\ref{t:Pet}. Let
$w\in W^P$ and consider quantum Schubert polynomial $C^P_w$
written as linear combination of $(\mathbf q,P)$-standard
monomials as in Section~\ref{s:Schub}. So
\begin{equation*}
C_{P,w}=\sum_{\Lambda\in\mathcal L_P}m_\Lambda
E_{P,\Lambda},\quad\qquad m_\Lambda\in\C.
\end{equation*}
In $E_{P,\Lambda}$ replace each factor $E^{(j)}_{P,i}$ with the
corresponding rational function $G^{n_j}_i$ to define $G_\Lambda$.
Then as function on $\mathcal Y_P$,
\begin{equation*}
\si^P_w=\sum_{\Lambda\in\mathcal L_P}m_\Lambda
G_{\Lambda}|_{\mathcal Y_P}.
\end{equation*}
We now use the `quantum straightening identity',
\cite{FoGePo:QSchub}{ Lemma 3.5},
\begin{equation*}\label{e:straightening}
E^{(j)}_{i} E^{(j)}_{l}=E^{(j+1)}_{i}E^{(j)}_{l}
-E^{(j)}_{i-1}E^{(j+1)}_{l+1}+E^{(j)}_{i-1}E^{(j)}_{l+1} +
q_j\left(E^{(j)}_{i-1} E^{(j-1)}_{l-1}- E^{(j-1)}_{i-2}
E^{(j)}_{l}\right)
\end{equation*}
to rewrite $\si^P_w$. Note that a factor $E^{(j)}_{i} E^{(j)}_{l}$
may occur in a $(\mathbf q,P)$-standard monomial $E_\Lambda$ only
if $j\in I^P$ and $j+1\notin I^P$. If we replace the
$E^{(j)}_i$'s by $G^{n_j}_i$ in the above identity and apply
Lemma~\ref{l:vanishing}, then we get
\begin{equation*}
\left(G^{n_j}_{i} G^{n_j}_{l}\right)|_{\Y_P}=
\left(G^{n_{j+1}}_{i}G^{n_j}_{l}
-G^{n_j}_{i-1}G^{n_{j+1}}_{l+1}+G^{n_j}_{i-1}G^{n_j}_{l+1}\right)|_{\Y_P}.
\end{equation*}

But the function $\sigma^B_w$ on $\Y_B$ (or equivalently the
rational function $\tilde\sigma_w\in \mathcal K(\Y)$) may be
obtained from the expression $\sum_{\Lambda\in\mathcal
L_P}m_\Lambda G_{\Lambda}$ we had for $\sigma^P_w$ by repeated
substitutions of the kind
\begin{equation}\label{e:substitutions}
G^{n_j}_{i} G^{n_j}_{l}\To G^{n_{j+1}}_{i}G^{n_j}_{l}
-G^{n_j}_{i-1}G^{n_{j+1}}_{l+1}+G^{n_j}_{i-1}G^{n_j}_{l+1},
\end{equation}
until the resulting expression has no more summands with factors
of type $G^{n_j}_{i} G^{n_j}_{l}$. (These transformations
correspond to the classical straightening identities which are
used to turn the $P$-standard monomial expansion of $c_w$ into
the $B$-standard monomial one). But the substitutions
\eqref{e:substitutions} do not affect the restriction to
$\mathcal Y_P$. So we are done.
\end{proof}

\begin{prop}\label{p:SchubClass}
For the Grassmannian permutation $w\in W^{P_m}$  define the
rational function $G_w$ on $G/B^-$ by
\begin{equation*}
G_w(gB^-):=\frac{(g\cdot v_{-\omega_m}\ |\ w\cdot v_{-\omega_m}) }
{(g\cdot v_{-\omega_m}\ |\ v_{-\omega_m})}
\end{equation*}
Then
\begin{equation*}
G_w|_{\Y}=\widetilde\si_w \ \in \mathcal K(\Y).
\end{equation*}
\end{prop}
\begin{proof}
By Proposition~\ref{p:SchubRestr} it suffices to show that
$G_w|_{\Y_B}$ coincides with $\sigma^B_w$. But this follows from
A. N. Kirillov's explicit formula for the corresponding quantum
Schubert polynomials, see Section~\ref{s:Kirillov}, together with
Peterson's Theorem~\ref{t:Pet}, and inspection of the matrix $u$
from \eqref{e:u} in the case where $P=B$.
\end{proof}

\begin{cor}\label{c:Grpos}

\begin{enumerate}
\item
If $y\in \Y_{B,>0}$, then for any $i\in I$ and $w\in W^{P_i}$ we
have $\si^B_w(y)>0$.
\item
If $y\in  \X_{P,>0}$ then $q^P_{i}(y)> 0$ for all $i=1\dotsc, k$.
\end{enumerate}
\end{cor}

\begin{proof}
(1) is an immediate corollary of Proposition~\ref{p:SchubClass},
since for any $g\in U^+$, $G_w(gB^-)$ is a quotient of nonzero
minors of $g$. Part (2) follows from
Proposition~\ref{p:SchubClass} along with Theorem~\ref{t:Pet}(3)
and Lemma~\ref{l:qpos}.
\end{proof}

\section{Proof of Theorem \ref{t:main}.(3)}\label{s:proof3}

We begin with a partial converse to Corollary~\ref{c:Grpos}.(1) in
the full flag variety case.
\begin{lem}\label{l:incl}
Let $y\in \Y_B$ such that  $\si^B_w(y)>0$ for all  $w\in W$. Then
$y\in \X_{B,>0}$. In other words,
$\Y_{B,>0}^{Schub}=\X_{B,>0}^{Schub}\subset \X_{B,>0}$. And
therefore also $X_{B,>0}^{Schub}\subset X_{B,>0}$.
\end{lem}

\begin{proof} By Lemma~\ref{l:qpos} we have that $q_i^B(y)>0$ for
all $i=1,\dotsc,k$. Therefore $y\in\X_B$. Now we may write $y=x
w_0 B^-$ for some $x\in X_B$. It remains to prove that $x\in
U^-_{>0}$. The positivity of all the quantum parameters $q^B_i$
implies by \eqref{e:Kos} that $\Delta_j(x)>0$ for all
$j=1,\dotsc,n-1$. Now by Proposition~\ref{p:SchubClass} the
positivity of the $\sigma^B_w$ for the Grassmannian permutations
$w$ of descent $d$ in $W$ implies the positivity of all the $d\x
d$ minors with column set $\{1,\dotsc, d\}$ and arbitrary row
sets. But this suffices to determine that $x$ is totally
positive, see e.g. \cite{BeFoZe:TotPos}.
\end{proof}

\begin{prop}\label{p:connected}
$X_{B,>0}^{Schub}=X_{B,>0}$.
\end{prop}

\begin{proof}
By  Lemma \ref{l:incl} we have the following commutative diagram
\begin{equation*}
\begin{matrix}X^{Schub}_{B,>0}&\hookrightarrow &X_{B,>0}\\
\qquad   \searrow & & \swarrow\qquad \\
      &\R_{>0}^{n-1}&
\end{matrix}
\end{equation*}
where the top row is clearly an open inclusion and the maps going
down are restrictions of $\pi^B$. By (1) of Theorem~\ref{t:main},
which is already proved, the left hand map to $\R_{>0}^{n-1}$ is a
homeomorphism. It follows from this and elementary point set
topology that $X^{Schub}_{B,>0}$ must be closed inside
$X_{B,>0}$. So it suffices to show that $X_{B,>0}$ is connected.

For an arbitrary element $u\in X$ and $t\in \R$, let
\begin{equation}\label{e:ut}
u_t:=\begin{pmatrix} 1& & & & \\
 ta_1& 1   &    &  &        \\
 t^2a_2   &t a_1 & \ddots&      &   \\
 \vdots   &     & \ddots & 1   &     \\
 t^{n-1}a_{n-1}   &\cdots     & t^2  a_2  &t a_1   & 1
\end{pmatrix}
\end{equation}
So $u_0=\Id$ and $u_1=u$, and if $u\in X_{B,>0}$, then so is $u_t$
for all positive $t$.

Let $u, u'\in X_{B,>0}$ be two arbitrary points. Consider the
paths
\begin{eqnarray*}
\gamma\ : {[0,1]\to X_{B,>0} ~,}& \gamma(t) =u u'_t\ \\
\gamma' : {[0,1]\to X_{B,>0} ~,}& \,\gamma'(t)=u_t u'.
\end{eqnarray*}
Note that these paths lie entirely in $X_{B,>0}$ since $X_{B,>0}$
is a semigroup (as the intersection of the group $X$ with the
semigroup $U^-_{>0}$). Since $\gamma$ and $\gamma'$ connect $u$
and $u'$, respectively, to $uu'$, it follows that $u$ and $u'$
lie in the same connected component of $X_{B,>0}$, and we are
done.
\end{proof}

\begin{cor}
$\Y_{P,>0}=\X_{P,>0}=\X_{P,>0}^{Schub}=\Y_{P,>0}^{Schub}$ and in
particular also $X_{P,>0}=X_{P,>0}^{Schub}$.
\end{cor}

\begin{proof}
The identity $\X_{P,>0}^{Schub}=\Y_{P,>0}^{Schub}$ follows from
Lemma~\ref{l:qpos}. It remains only to show that
$X_{P,>0}=X_{P,>0}^{Schub}$. We begin with the inclusion
$\subseteq$. Let $X_{\ge 0}=X\cap U^-_{\ge 0}$. Then clearly
\begin{equation}\label{e:closures}
\overline {X_{B,> 0}}\subseteq X_{\ge 0}
\end{equation}
is an inclusion of closed subsemigroups of $U^-$. We show that
this is actually an equality. Suppose $x\in X_{\ge 0}$, then for
any $u\in X_{B,>0}$ and $u_t$ defined as in \eqref{e:ut}, the
curve $t\mapsto x(t)=xu_t$ starts at $x(0)=x$ and lies in
$X_{B,>0}$ for all $t>0$. Therefore $x\in \overline {X_{B,>0}}$
as desired. As a consequence, using
Proposition~\ref{p:connected}, we have
\begin{equation}\label{e:closures2}
X_{P,>0}=X_{P}\cap \overline{X_{B,>0}}=X_{P}\cap
\overline{X^{Schub}_{B,>0}}.
\end{equation}
Now consider the Schubert classes $\si^P_{w}\in qH^*(G/P)$ as
functions on $X_P$. By Proposition \ref{p:SchubRestr},
$\si^P_{w}=\tilde\si_{w}|_{X_P}$, and $\widetilde\si_{w}$ takes
positive values on $X^{Schub}_{B,>0}$. Let us choose $x\in
X_{P,>0}$. Then by \eqref{e:closures2} we have also $\si^P_w(x)\ge
0$ for all $w\in W^P$. On the other hand
$Q:=\pi^P(x)=(q^P_1(x),\dotsc,q^P_k(x))\in \R_{>0}^k$ by
Corollary~\ref{c:Grpos}. But we have seen in
Section~\ref{s:proof1} that there is only one Schubert nonnegative
point in the fiber $(\pi^P)\inv(Q)$, and that that one is strictly
positive. Thus in fact $\si^P_w(x)> 0$ for all $w\in W^P$ and
$X_{P,>0}\subset X^{Schub}_{P,>0}$.

It remains to show that $X_{P,>0}\hookrightarrow X_{P,>0}^{Schub}$
is surjective. Consider again the proper map
\begin{equation*}
\Delta=(\Delta_1,\dotsc,\Delta_{n-1}): X\to \C^{n-1}
\end{equation*}
defined in Section~\ref{s:Toeplitz}. Its restriction $\Delta_{\ge
0}=(\Delta_1,\dotsc,\Delta_{n-1}): X_{\ge 0}\to (\R_{\ge
0})^{n-1}$ is surjective, since the image must be closed and
contain $\Delta_{\ge 0}(X_{B,>0})=\R_{>0}^{n-1}$.

From Theorem~\ref{t:main}(1) along with Lemma~\ref{l:GenKos} we
know that the further `restriction' of $\Delta$,
\begin{equation*}
\Delta^P_{>0}=(\Delta_{n_1},\dotsc,\Delta_{n_k}):
X^{Schub}_{P,>0}\to (\R_{>0})^k,
\end{equation*}
is bijective. So we have the following diagram,
\begin{equation*}
\begin{matrix}X_{P,>0}&\hookrightarrow &X^{Schub}_{P,>0}\\
       \qquad \searrow & & \swarrow\sim\qquad \\
      &\R_{>0}^{k}&
\end{matrix}
\end{equation*}
where the downward arrows are given by $\Delta^P_{>0}$ and its
restriction. By the surjectivity of $\Delta_{\ge 0}$ we also have
that the left hand map is surjective. This implies the desired
equality, $X_{P,>0}= X^{Schub}_{P,>0}$.
\end{proof}

\def\cprime{$'$}
\providecommand{\bysame}{\leavevmode\hbox
to3em{\hrulefill}\thinspace}


\end{document}